\documentclass[a4paper,american,reqno]{amsart}
\usepackage{babel}
\usepackage[utf8]{inputenc}
\usepackage[T1]{fontenc}
\usepackage{microtype}
\usepackage{enumitem}
\usepackage{booktabs}
\usepackage{pgfplots}
\usepackage{relsize}
\usepackage{tikz}
\usepackage{xspace}
\usepackage{csquotes}
\usepackage[colorlinks,linkcolor=black]{hyperref} 
\hypersetup{
 colorlinks, citecolor=black, filecolor=black, linkcolor=black, urlcolor=black
}
\usepackage{algorithm}
\usepackage{algorithmic} 
\usepackage{float}
\usepackage{array}
\usepackage{amssymb}
\usepackage{pdfpages}
\usepackage{setspace}
\usepackage{filecontents}

\newtheorem{lemma}{Lemma}
\newtheorem{proposition}{Proposition}
\newtheorem{theorem}{Theorem}

\newtheorem{assumption}{Assumption}

\theoremstyle{definition}

\newtheorem{remark}{Remark}
\newtheorem{example}{Example}

\makeatletter
\patchcmd{\@settitle}{\uppercasenonmath\@title}{\scshape\large}{}{}
\patchcmd{\@setauthors}{\MakeUppercase}{\scshape\normalsize}{}{}
\makeatother

\vfuzz \hfuzz
\renewcommand{\:}{\colon}
\renewcommand{\j}{\jmath}

\newcommand{\RR}{\mathbb{R}}
\newcommand{\NN}{\mathbb{N}}

\newcommand{\SG}{S}
\newcommand{\tf}{t_{\mathrm{f}}}
\newcommand{\Ut}{U_{[t_0,\tf]}}
\newcommand{\Sigmat}{\Sigma_{[t_0,\tf]}}
\newcommand{\Wt}{W_{[t_0,\tf]}}
\newcommand{\CM}{\ensuremath{Q}}
\newcommand{\CT}{\ensuremath{\mathcal{T}}}
\newcommand{\dual}{\ensuremath{{Y^{{}^*}\!\!,\, Y}}}

\DeclareMathOperator{\cvto}{\curvearrowright}

\DeclareMathOperator{\laplace}{\Delta}
\DeclareMathOperator{\Vp}{\raisebox{-0.5em}{\ensuremath{+}}\kern-1.0ex\bigvee}
\DeclareMathOperator{\Vm}{\raisebox{-0.5em}{\ensuremath{--}}\kern-1.0ex\bigvee}

\newcommand{\eg}{e.g.}
\newcommand{\ie}{i.e.}
\newcommand{\cf}{cf.}

\begin{document}

\title[Stability and Optimal Control of Switching PDE-Dynamical Systems]{Stability and Optimal Control of\\ Switching PDE-Dynamical Systems}
\author[F.\,M. Hante]{Falk~M.~Hante$^{1}$}
\address{$^1$Friedrich-Alexander-Universität Erlangen-Nürnberg,
  Lehrstuhl für Angewandte Mathematik 2,
  Cauerstr.~11,
  91058~Erlangen,
  Germany.}
\date{February 15, 2018. An earlier and extended version of this manuscript was submitted by the author as \emph{Habilitationsschrift} to the Faculty of Sciences at the 
University of Erlangen-Nürnberg, Germany, on Mai 15, 2017.}

\begin{abstract}Selected results for the stability and optimal control of abstract switched systems in Banach and Hilbert space are reviewed. The dynamics are typically given in a piecewise sense by a family of nonlinearly perturbed evolutions of strongly continuous semigroups. Stability refers to characterizations of asymptotic decay of solutions that holds uniformly for certain classes of switching signals for time going to infinity. Optimal control refers to the minimization of costs associated to solutions by appropriately selecting switching signals. Selected numerical results verify and visualize some of the available theory. \end{abstract}

\maketitle

\section{Introduction} \label{sec:intro}
This manuscript summarizes selected contributions concerning abstract evolution systems of the form
\begin{equation}\label{eq:sysGeneral}
\begin{aligned}
 \frac{d}{dt}y &= A^{\sigma(t)} y + f^{\sigma(t)}(y)~\text{on}~Y,\quad~t>0,\\
 y(0)&=y_0,
\end{aligned}
\end{equation}
where for some index set $Q$ of finite or infinite cardinality $\sigma(\cdot)\: [0,\infty) \to Q$ is a switching signal conducting
the temporal evolution initiated in $y_0 \in Y$ among families of possibly unbounded linear 
operators $\{A^j\: D(A^j) \subset Y \to Y\}_{j \in Q}$ and nonlinear perturbations $\{f^j\: Y \to Y\}_{j \in Q}$.

Unless explicitly stated, it is assumed that $Y$ is a reflexive Banach space, $\sigma$ is piecewise constant and that, for each $j \in Q$, 
the operator $A^j$ is the generator of a strongly continuous semigroup of bounded linear operators $\SG^j(t)$ on $Y$ with a perturbation
$f^j$ being Lipschitz continuous and appropriately bounded. For precise technical assumptions concerning, \eg, the type of bound on 
$f^j$ it is referred to the respective original works. 

The switching signal $\sigma(\cdot)$ in \eqref{eq:sysGeneral} can be identified with a sequence of switching times
\begin{equation}
 0=\tau_0 \leq \tau_1 \leq \tau_2 \leq \tau_3 \leq \ldots \quad\text{in}~[0,\infty)
\end{equation}
and a sequence of modes 
\begin{equation}
 j_0,~j_1,~j_2,~j_3,~\ldots\quad \text{in}~Q
\end{equation}
using the relation
\begin{equation}
 \sigma(t)=j_{k},\quad t\in [\tau_{k},\tau_{k+1}),~k \in \NN_0.
\end{equation}
A solution of the abstract evolution equation \eqref{eq:sysGeneral} is then to be understood 
in the sense of mild solutions given by a continuous function $y \in C([0,\infty);Y)$ satisfying the initial 
condition $y(0)=y_0$ and the variation of constants formula in the piecewise sense
\begin{equation}\label{eq:varconst}
y(t) = \SG^{j_{k}}(t-\tau_k)y(\tau_{k}) + \int_{\tau_{k}}^{t} \SG^{j_{k}}(t-\tau_k-s)f^{j_{k}}(y(s))\,ds,~t \in [\tau_{k},\tau_{k+1}),~k \in \NN_0.
\end{equation}

According to this definition, the switching system \eqref{eq:sysGeneral} covers in an abstract sense many evolution problems involving ordinary differential equations (ODEs), linear delay differential equations (DDEs), and linear partial differential equations (PDEs). We note that in cases of PDEs, the switching signal $\sigma$ may also switch the principle part of the equation. This situation is explicitly considered in the publications \cite{HanteSigalotti2011}, \cite{AminHanteBayen2012}, and \cite{RuefflerHante2016} and is otherwise rarely addressed in the available literature concerning switched systems. 
The publications \cite{HanteSigalottiTucsnak2012}, \cite{HanteSager2013}, \cite{Hante2017}, and \cite{GugatHante2017} address the special case that $A^{\sigma(\cdot)}$ is $\sigma$-invariant, \ie, the switching signal $\sigma$ solely acts on the perturbation. In these cases, we write $A^{\sigma(\cdot)}=A$ and assume that $A$ is the infinitesimal generator of a strongly continuous semigroup $\SG(t)$ on $Y$. Moreover, in these cases, \eqref{eq:varconst} simplifies to
\begin{equation}\label{eq:varconst2}
y(t) = \SG(t)y(0) + \int_{0}^{t} \SG(t-s)f^{\sigma(t)}(y(s))\,ds,\quad t>0,
\end{equation}
and we may also extend the class of switching signals $\sigma(\cdot)$ to measurable functions $\sigma\: [0,\infty) \to Q$.

Necessary and sufficient conditions on a linear operator $A$ to be the infinitesimal generator of a strongly continuous semigroup go back to E. Hille, 
G. Lumer, R.\,S. Phillips, K. Yosida, and M.\,H. Stone and can be found in many textbooks; see, \eg, \cite{Balakrishnan1981,BensoussanDaPratoDelfourMitter1992,EngelNagel1999,LiYong1995,Pazy1983}. 
To illustrate the range of applicability we discuss some academic examples of switching PDE-dynamical systems being covered by the abstract 
setting~\eqref{eq:sysGeneral}.

\begin{example}\label{ex:schoedinger}
Let $N \geq 1$, $Q=\{0,1\}$, and consider over the complex field an $N$-dimensional Schrödinger equation with on/off damping
\begin{equation}\label{eq:schroedinger-N}
\begin{aligned}
 	i\,y_t(t,x) +\laplace y(t,x) + i\,\sigma(t)\,d(x)^2 y(t,x)& = 0, ~&(t,x)\in(0,\infty)\times \Omega,\\
 	y(t,x)&=0	,					 ~& t\in (0,\infty)\times \partial \Omega,\\
 	y(0,x) &= y_0(x),					 ~& t\in \Omega,
\end{aligned}
\end{equation}
where $\Omega$ is a bounded domain in $\RR^N$ with regular boundary $\partial\Omega$ and $d \in L^\infty(\Omega)$.
Let $A$ be the linear operator on $Y=L^2(\Omega)$ defined by
\begin{equation}
 Ay= i\,\Delta y\quad~\text{on}~D(A)=H^2(\Omega) \cap H^1_0(\Omega),
\end{equation}
and let $B\: Y \to Y$ be the multiplication operator $By=d\,y$.
It is well known that $D(A)$ is dense in $Y$ and easy to verify that $A$ is skew-adjoint, \ie, $A^*=-A$. Hence,
it follows from the theorem of Stone that $A$ is the generator of a strongly continuous group of
unitary operators on $Y$; see, \eg, \cite{Rudin1991}.
Then, setting  
\begin{equation}\label{eq:shroedingerabstract}
\begin{aligned}
A^j=A,\quad D(A^j)=D(A), \quad f^j(y)=j\,BB^*y,\quad j \in Q,
\end{aligned}
\end{equation}
where $B^*$ denotes the adjoint of $B$, problem \eqref{eq:schroedinger-N} is equivalent to \eqref{eq:sysGeneral}. \qed
\end{example}

\begin{example}\label{ex:heat} 
Let $N \geq 1$, $Q=\{0,1\}$, and consider over the field of real numbers an $N$-dimensional heat equation on a 
bounded domain $\Omega$ with a regular boundary $\partial\Omega$ being controlled by applying a lumped 
control $u$ on alternating control subdomains~$\omega_j \subset \Omega$
\begin{equation}\label{eq:heat-N}
\begin{aligned}
 	y_t(t,x) -\laplace y(t,x) & = \chi_{\omega_{\sigma(t)}}(x)u(t), ~&(t,x)\in(0,\infty)\times \Omega,\\
 	y(t,x)&=0	,					 ~& t\in (0,\infty)\times \partial \Omega,\\
 	y(0,x) &= y_0(x),					 ~& t\in \Omega,
\end{aligned}
\end{equation}
where $\chi_{\omega_{j}}$ denotes the characteristic function of $\omega_j$. Let $A$ be the linear operator on $Y=L^2(\Omega)$ defined by the Dirichlet--Laplace operator
\begin{equation}
 Ay= \Delta y\quad~\text{on}~D(A)=H^2(\Omega) \cap H^1_0(\Omega).
\end{equation}
One can use the Hille--Yosida theorem to conclude that $A$ generates a strongly continuous semigroup. Further, the generated 
semigroup $\SG(t)$ can be shown to be contractive and analytic; see, \eg, \cite{Pazy1983}.
Then, setting
\begin{equation}\label{eq:heatabstract}
\begin{aligned}
A^j=A,\quad D(A^j)=D(A), \quad f^j(y)=f^j(y,u)=\chi_{\omega_{j}}u,\quad j \in Q,
\end{aligned}
\end{equation}
problem \eqref{eq:heat-N} is equivalent to \eqref{eq:sysGeneral}. \qed
\end{example}

\begin{example}\label{ex:wave}
Let $N \geq 1$, $Q=\{0,1\}$, and consider over the field of real numbers an $N$-dimensional wave equation with on/off damping
\begin{equation}\label{eq:wellengleichung-N} 
\begin{aligned}
v_{tt}(t,x)&= \Delta v(t,x)-\sigma(t)\,d^2(x) v_t(t,x),\quad &(t,x)\in (0,\infty) \times \Omega,\\
v(t,x)& = 0,\quad &(t,x)\in (0,\infty) \times  \partial \Omega,\\
v(0,x)&= y_0(x),\quad v_t(0,x) = y_1(x), &x\in \Omega,\\
\end{aligned}
\end{equation}
where $\Omega$ is a bounded domain in $\RR^N$ with regular boundary $\partial\Omega$ and $d \in L^\infty(\Omega)$.
Let $A$ be the linear operator on $Y=H^1(\Omega) \times L^2(\Omega)$ defined by
\begin{equation}
\begin{aligned}
 D(A) &= \left\{ \begin{pmatrix}y_1\\y_2\end{pmatrix} : y_1 \in H^2(\Omega)\cap H^1_0(\Omega),~y_2 \in H^1(\Omega)\right\},\\
 A\begin{pmatrix}y_1\\y_2\end{pmatrix} &= \begin{pmatrix}y_2\\ \Delta y_1\end{pmatrix}\quad\text{on}~D(A).
\end{aligned}
\end{equation}
It is well-known that $D(A)$ is dense in $Y$, $A$ is closed, and that $A$ is maximally dissipative for $Y$ equipped with the norm
\begin{equation}
\left\|\begin{pmatrix}y_1\\y_2\end{pmatrix}\right\|_Y^2=\|\nabla y_1\|^2_{L^2(\Omega)}+\| y_2\|^2_{L^2(\Omega)},
\end{equation}
see, \eg, \cite{TucsnakWeiss2009}. Hence, the Lumer--Phillips theorem yields that $A$ is the generator of a strongly 
continuous semigroup of contractions on $Y$. Further, let $B$ be the linear operator on $Y$ defined by 
\begin{equation}
B\begin{pmatrix}y_1\\y_2\end{pmatrix}(x)=\begin{pmatrix}0\\d(x)\,y_2(x)\end{pmatrix},\quad x\in \Omega. 
\end{equation}
Then, setting  
\begin{equation}\label{eq:waveabstract}
\begin{aligned}
A^j=A,\quad D(A^j)=D(A), \quad f^j(y)=j\,BB^*y,\quad j \in Q,
\end{aligned}
\end{equation}
where $B^*$ denotes the adjoint of $B$, problem \eqref{eq:wellengleichung-N} 
is equivalent to \eqref{eq:sysGeneral}. \qed
\end{example}

\begin{example}\label{ex:transport}
Let $n \geq 1$, $Q=\{1,\ldots,M\}$, and consider over the field of real numbers an $n$-dimensional system of transport equations on some interval $[a,b]$ 
with switching transport velocities and switching boundary conditions
\begin{equation}\label{eq:transport-n}
\begin{aligned}
 &y_t(t,x) + \Lambda^{\sigma(t)}(x) y_x(t,x) = f^{\sigma(t)}(y(t,x)),\quad &(t,x)\in(0,\infty) \times (a,b),\\
 &y_{II}(t,a)=G_L^{\sigma(t)} y_{I}(t,a),~y_{I}(t,b)=G_R^{\sigma(t)} y_{II}(t,b),\quad &t \in (0,\infty),\\
 &y(0,x)=y_0(x),\quad &x \in (a,b),
\end{aligned}
\end{equation}
where on $[a,b]$, for all $j \in Q$ and some $m \in \{1,\ldots,n\}$, $\Lambda^j(x)$ are diagonal and continuously differentiable matrix functions 
$\Lambda^j(x) = \mathrm{diag}\,(\lambda_1^j(x),\ldots,\lambda_n^j(x))$ satisfying
 \begin{equation}\label{eq:stricthyperb}
  \lambda_1^j(x)<\ldots<\lambda_{m}^j(x)<0<\lambda_{m+1}^j(x)<\ldots<\lambda_n^j(x),\quad x \in [a,b],
 \end{equation}
where the state vector partitions as $y(t,x)=(y_I(t,x)^\top,y_{II}(t,x)^\top)^\top$ with
\begin{equation}
 y_I(t,x)=(y_1(t,x),\ldots,y_m(t,x))^\top,\quad y_{II}(t,x)=(y_{m+1}(t,x),\ldots,y_n(t,x))^\top,
\end{equation}
and where $G_L^j$ and $G_R^j$ are matrices of dimensions ${(n-m)\times m}$ and ${m\times (n-m)}$, respectively. 
For all $j \in Q$, let $A^j$ be the linear operator on $Y=L^2(a,b;\RR^n)$ defined by
\begin{equation}
\begin{aligned}
 &D(A^j)=\biggl\{\begin{pmatrix}y_I\\y_{II}\end{pmatrix} \in H^1(a,b;\RR^m) \times H^1(a,b;\RR^{(n-m)}) : \\
 &\qquad\qquad\qquad y_{II}(a)=G_L^j y_{I}(a),~y_{I}(b)=G_R^j y_{II}(b)\biggr\},\\ 
 &A^j y(x) = - \Lambda^j(x) \frac{d}{dx}y,\quad x \in (a,b),~y \in D(A^j).
\end{aligned}
\end{equation}
In this case, one can use the method of characteristics to verify that for each $j \in Q$, $A^j$ generates on $Y$ a strongly continuous 
semigroup \cite{Russell1978}. Hence, problem \eqref{eq:transport-n} is equivalent to \eqref{eq:sysGeneral}. \qed
\end{example}

\begin{figure}
\includegraphics[width=0.8\textwidth]{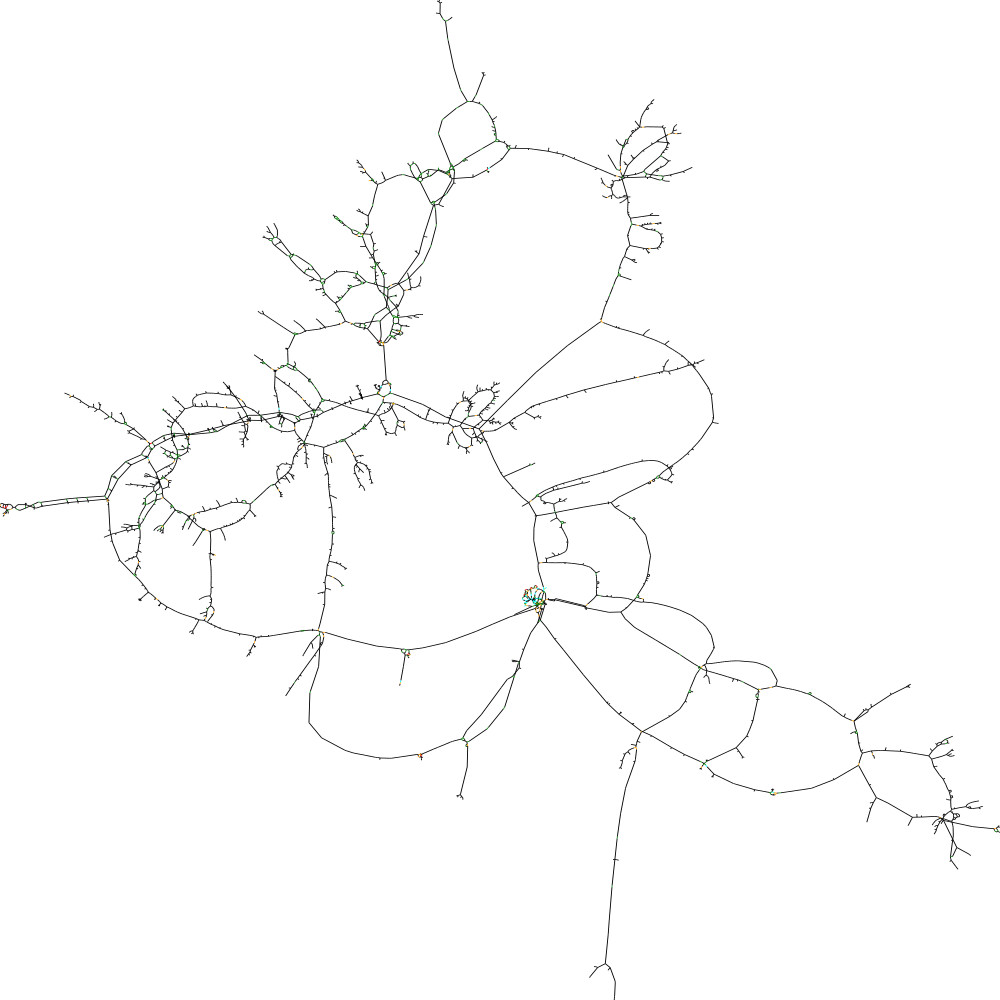}
\caption{The gas network example GasLib-4197 \cite{GasLib}. The network contains $426$ valves that can be opened and closed to control the gas flow. In subsonic regimes the pressure and flow 
distribution over time can be modeled by  semilinear Euler gas equations using a switching system as in Example~\ref{ex:transport}. \hfill {\footnotesize Source: \url{http://gaslib.zib.de/}}}
\label{fig:gas} 
\end{figure}

We note that the hypothesis for the examples above can be met for many physical systems; see, \eg, \cite{CourantHilbert1962}.
Example~\ref{ex:heat}, for instance, is motivated by linearized problems in thermal manufacturing systems where the discrete mode~$j \in Q$
models actuator placements \cite{IftimeDemetriou2009} and the assumed diagonal form of \eqref{eq:transport-n} in Example~\ref{ex:transport} 
can be obtained for semilinear models considered for the flow of gas or fresh water in pipe networks as well as drainage or sewer systems 
in open canals, where the discrete mode $j \in Q$ corresponds to some fixed valve position subject to switching \cite{HanteEtAl2017}; 
\cf~Figure~\ref{fig:gas}. 

The contributions summarized below concern two in some sense complementary problems for switching systems as those above. The first problem
is related to the asymptotic stability of solutions of \eqref{eq:sysGeneral}, where the main focus lies on conditions that guarantee the respective 
asymptotic decay properties uniformly for certain subclasses of switching signals $\sigma(\cdot)$. Available results in this direction for 
finite-dimensional switched linear systems are surveyed in \cite{ShortenEtAl2007}. In particular, it is well-known that stability conditions for 
all subsystems with $\sigma(t) \equiv j$ for some $j  \in Q$ is in general not sufficient to ensure uniform stability of the switched system.
Section~\ref{sec:stability} summarizes the contributions of the author made in \cite{HanteSigalotti2011}, \cite{AminHanteBayen2012}, and \cite{HanteSigalottiTucsnak2012} 
for certain infinite-dimensional linear switched systems as subclasses of~\eqref{eq:sysGeneral} concerning uniform stability conditions.
The second problem is related to minimization of costs associated to the solution of the switched systems \eqref{eq:sysGeneral} by appropriately 
selecting the switching signals $\sigma(\cdot)$ possibly in addition to further continuous valued control functions $u$. Available results in this 
direction for ODE-dynamical problems are surveyed in \cite{LinAntsaklis2009,Sager2009}. In particular, it is well-known that certain stabilization 
problems may be posed as optimal control problems and that in this context, stability conditions for all subsystems 
with $\sigma(t) \equiv j$ for some $j  \in Q$ are in general not necessary to obtain asymptotic stability for some switching signal $\sigma(\cdot)$.
Section~\ref{sec:optimalswitching} summarizes the contributions of the author made in \cite{HanteSager2013}, \cite{Hante2017}, \cite{GugatHante2017}, 
and \cite{RuefflerHante2016} for PDE-dynamical problems in the framework of~\eqref{eq:sysGeneral} concerning such switching optimal control problems. In order to verify and visualize 
the theoretical contributions, selected numerical results from \cite{AminHanteBayen2012} concerning the stability of switched systems and 
from \cite{Hante2017} concerning optimal switching control are summarized at the end of Section~\ref{sec:stability} and Section~\ref{sec:optimalswitching}, respectively. A conclusion is drawn in Section~\ref{sec:conclusion}.

\section{Stability of switching PDE-dynamical systems}\label{sec:stability}
Consider a switched system \eqref{eq:sysGeneral} in the homogeneous case, \ie, a system of the form
\begin{equation}\label{eq:ssOpForm}
\begin{aligned}
 \frac{d}{dt}\,y(t)	&= A^{\sigma(t)} y(t),\quad t > 0,\\
 	y(0)		&= y_0 \in Y.
\end{aligned}
\end{equation}
In this case, the origin $y=0$ is an equilibrium solution and we are interested in finding conditions that ensure the convergence 
of the solution $y(t)$ to the origin as $t$ tends to infinity uniformly for switching signals $\sigma(\cdot)$ in a particular class of 
signals. It is well-known that even in the finite-dimensional case stability conditions for each individual $A^j$, $j \in Q$, are not 
sufficient for such a stability property to hold, \eg, for all possible signals $\sigma(\cdot)$, see \cite{ShortenEtAl2007}.

In Section~\ref{sec:GUES} we review results from \cite{HanteSigalotti2011} and \cite{AminHanteBayen2012} concerning the general system~\eqref{eq:ssOpForm} being global uniformly exponentially stable, \ie, stability being asymptotically with an exponential rate of convergence
and uniformly with respect to all possible signals and all possible initial conditions.
In Section~\ref{sec:stabother} we summarize results from \cite{HanteSigalottiTucsnak2012} concerning exponential, asymptotic, and weak
asymptotic stability of a system \eqref{eq:ssOpForm} for the special case that $Y$ is a Hilbert space and $A^{j}=A-jBB^*$ with a 
dissipative operator $A$ and a bounded operator $B$ and uniformity with respect to signals satisfying persistent excitation 
conditions. The motivation behind this particular setting is in the interpretation of the switching as intermittencies of classical
feedback control for dissipative systems.

\subsection{Global uniform exponential stability of switched systems}\label{sec:GUES}
For finite-dimensional switched systems, \ie, where $A^j$, $j \in Q$, are matrices and 
\begin{equation}
S^j(t)=\sum_{k=0}^\infty \frac{t^k\left(A^j\right)^k}{k!}
\end{equation}
is the matrix exponential, it is well-known that a necessary and sufficient condition for global uniform exponential stability 
is the existence of a common Lyapunov function \cite{MolchanovPyatnitskiy1989,MolchanovPyatnitskiy1986}. It is even necessary and 
sufficient for the global uniform asymptotic stability \cite{LinSontagWang1996} in case of nonlinear subsystems.
Furthermore, in finite dimension and if the switched system has finitely many modes, it is known that the common Lyapunov function can be 
taken polyhedral or polynomial; see \cite{Blanchini1994,BlanchiniMiani1999,DayawansaMartin1999} and also \cite{BraytonTong1980}
for a discrete-time version. A special role in the switched control literature has been played by common quadratic Lyapunov functions because 
their existence can be tested rather efficiently. It is known, however, that the existence of a common quadratic Lyapunov function 
is not necessary for the global uniform exponential stability of a linear switched system with finitely many modes. Moreover, there exists 
no uniform upper bound on the minimal degree of a common polynomial Lyapunov function \cite{MasonBoscainChitour2006}. 
It is referred to \cite{LinAntsaklis2009,ShortenEtAl2007} for surveys the available results. The characterization of exponential stability 
for a single linear dynamical system on infinite-dimensional Banach and Hilbert spaces dates back to Datko \cite{Datko1968} and Pazy \cite{Pazy1972} 
and has, since then, seen a broad range of applications in control theory for partial differential equations; see, \eg, \cite{TucsnakWeiss2009}.

In \cite{HanteSigalotti2011}, some of the above known results in finite dimension have been extended to in\-fi\-ni\-te-di\-men\-sio\-nal 
switched systems of the type \eqref{eq:ssOpForm} as a homogeneous switched system in the framework of \eqref{eq:sysGeneral}. The main result that 
has been obtained extends the necessary and sufficient condition for global uniform exponential stability in terms of existence of a common Lyapunov 
function to a general Banach space setting. The result can be summarized as follows.
\begin{theorem}\label{thm0:GUES} The following three conditions are equivalent:
\begin{itemize}
 \item[(A)] There exist two constants $K \geq 1$ and $\mu>0$ such that  
for every $\sigma(\cdot)$ and every $y_0$ the solution $y(\cdot)$ to \eqref{eq:ssOpForm} satisfies
\begin{equation}
 \|y(t)\|_{Y} \leq Ke^{-\mu t}\|y_0\|_{Y}, \quad t \geq 0.
\end{equation}
\item[(B)] There exist two constants $M \geq 1$ and $\omega>0$
such that, for every $\sigma(\cdot)$ and every $y_0$,  the solution $y(\cdot)$ to \eqref{eq:ssOpForm} satisfies
\begin{equation}\label{thm0:UEBound}
 \|y(t)\|_{Y} \leq Me^{\omega t}\|y_0\|_{Y}, \quad t \geq 0,
\end{equation}
and there exists $V\: Y \to [0,\infty)$
such that $\sqrt{V(\cdot)}$ is a norm on $Y$, 
	\begin{equation}\label{thm0:VBound}
 		V(y) \leq C \|y\|_Y^2, \quad y\in Y,
	\end{equation}
	for a constant $C>0$, and
	\begin{equation}\label{thm0:LjVest}
 		\liminf_{t \downarrow 0} \frac{V(\SG^j(t)y) - V(y)}{t} \leq -\|y\|_Y^2,\quad
		j \in Q,~y\in Y.
	\end{equation}
\item[(C)] There exists $V\: Y \to [0,\infty)$
such that $\sqrt{V(\cdot)}$ is a norm on $Y$, 
	\begin{equation}\label{thm00:VBound}
 		c\|y\|_Y^2\leq V(y) \leq C \|y\|_Y^2, \quad y\in Y,
	\end{equation}
	for some constants $c,C>0$, and
	\begin{equation}\label{thm00:LjVest}
 		\liminf_{t \downarrow 0} \frac{V(\SG^j(t)y) - V(y)}{t} \leq -\|y\|_Y^2,\quad j \in Q,~y\in Y.
	\end{equation}
\end{itemize}
\end{theorem}

The construction of a common Lyapunov function used in \cite{HanteSigalotti2011} satisfying (B) under the assumption that (A) holds true uses the  
candidate function
\begin{equation*}
V(y_0)=\sup\left\{ \int_0^\infty \|y(t)\|^2 dt : y(\cdot)~\text{solution of \eqref{eq:ssOpForm} for some}~\sigma\right\}.
\end{equation*}
Alternatively, one can take $V(y_0)=\int_0^\infty \sup_{\sigma(\cdot)}\|y(t)\|^2 dt$, as done in \cite{HanteSigalotti2010}.
The construction of a Lyapunov function satisfying (C) under the assumption that (A) holds true is similar, if 
\eqref{eq:ssOpForm} is augmented with a further mode $\SG^{j^*}(t)=e^{-\mu t}I$, where $\mu>0$ is the constant appearing in (A) 
and $I$ denotes the identity on $Y$, and to consider all the solutions to this augmented system in the definition of $V$.

\begin{remark}
We note that the equivalence between (A) and (C) extends to infinite-dimensional systems the result of \cite{MolchanovPyatnitskiy1986},
and that it shows that the assumption of compactness of $\{A^j\mid j\in Q\}$ that is typically made in the finite-dimensional setting 
is not needed. The conditions \eqref{thm0:VBound} and \eqref{thm00:VBound} are redundant in the case of finite-dimensional systems, because 
$\sqrt{V(\cdot)}$ and $\|\cdot\|_Y$ are comparable by compactness of the unit sphere. Hence, condition \eqref{thm0:UEBound} in (B) could 
be dropped for finite-dimensional systems. This is not the case for infinite-dimensional ones, as illustrated in \cite{HanteSigalotti2011} 
by an example. We further note that in the case of an exponentially stable single mode ($Q=\{0\}$), it was observed by Pazy \cite{Pazy1972}
that $y \mapsto \int_0^\infty \|\SG_0(t)y\|^2\,dt$ defines a Lyapunov function that is comparable with the squared norm if and only if $\SG_0$ 
extends to an exponentially stable strongly continuous group. The above result shows that in this context, as a consequence of the 
implication (A)~$\Rightarrow$~(C), even if $\SG_0$ does not admit an extension to a group, a Lyapunov function comparable with the 
squared norm can still be found.
\end{remark}

From the point of view of applications, condition (B), imposing less conditions on $V$ than (C), is better suited for establishing that (A) holds 
(although the uniform exponential growth boundedness needs also be proved). On the other hand, the implication (A)~$\Rightarrow$~(C) can be used to 
select a Lyapunov function with tighter requirements. 

Concerning the regularity of the Lyapunov functions obtained through the construction in \cite{HanteSigalotti2011}, they are always convex and continuous 
because $\sqrt{V(\cdot)}$ is a norm. In the special case in which $Y$ is a separable Hilbert space, we obtained in \cite{HanteSigalotti2011} also the 
following additional regularity properties.

\begin{proposition}\label{prop:HilbertGUES}
Let $Y$ be a separable Hilbert space and assume that (A) in Theorem~\ref{thm0:GUES} holds. Then (C) in Theorem~\ref{thm0:GUES} holds with $V(\cdot)$ being a 
directionally Fr\'{e}chet differentiable function.
\end{proposition}

While the above theory is rather general, more specific conditions can be found guaranteeing global uniform exponential stability for particular systems.
We demonstrate this for hyperbolic initial boundary value problems from Example~\ref{ex:transport}.

\begin{theorem}\label{thm:GUEShyperbolic} Consider problem \eqref{eq:transport-n} with $f^j\equiv 0$ and assume that the matrices $G_L^j$ and $G_R^j$ satisfy
the spectral radius condition
\begin{equation}
 \max_{j,j'\in Q}~\rho\left(\begin{bmatrix} 0 & |G_R^{j'}| \\ |G_L^j| & 0 \end{bmatrix}\right)<1.
\end{equation}
Then, for any $y_0 \in L^\infty(a,b;\RR^n)$, the solution $y(\cdot)$ of \eqref{eq:ssOpForm} corresponding to the solution of \eqref{eq:transport-n}
satisfies
\begin{equation}\label{eq:GUEShyper}
 \|y(t)\|_{Y} \leq Ke^{-\mu t}\|y_0\|_{Y}, \quad t \geq 0,
\end{equation}
for some constants $K\geq 1$ and $\mu > 0$. 
\end{theorem}
The proof in \cite{AminHanteBayen2012} is based on recursively estimating the norm of the solution using the method of characteristics. The 
spectral radius condition thereby ensures that all possible combinations of boundary reflections are nonamplifying in the classical sense 
of Tatsien Li \cite{TTLi85}.
\begin{remark} The implication (A)~$\Rightarrow$~(C) of Theorem~\ref{thm0:GUES} yields that under the conditions of Theorem~\ref{thm:GUEShyperbolic},
there exists a common Lyapunov function for the problem in Example~\ref{ex:transport}. However, for $|Q|>1$ a closed form expression of such 
a Lyapunov function for example similar to those utilized in the single mode case in \cite{Coron07} is not yet known.
\end{remark}
\begin{remark}\label{rem:GUEShyperbolicInhomo}
The analysis in \cite{AminHanteBayen2012} also addresses the global uniform stability under similar conditions as those in Theorem~\ref{thm:GUEShyperbolic}
for well-posed non-diagonal inhomogeneous systems of the form
\begin{equation}
 A^j  y (x) = - \Lambda^j(x) \frac{d}{dx}y + B(x)y,\quad x \in (a,b),~y \in D(A^j),
\end{equation}
with $\Lambda^j \in \RR^{n\times n}$ under a commutativity assumption
\begin{equation}\label{eq:commute}
\Lambda^j(x)\Lambda^{j'}(x)=\Lambda^{j'}(x)\Lambda^{j}(x),\quad x\in[a,b],~j,j'\in Q,
\end{equation}
and a smallness condition on $\|B(\cdot)\|_\infty$. 
\end{remark}

A numerical example illustrating the results from Theorem~\ref{thm:GUEShyperbolic} and Remark~\ref{rem:GUEShyperbolicInhomo} is discussed in Section~\ref{sec:numresultsstab}.

\subsection{Stability of dissipative systems with intermittent damping}\label{sec:stabother}The notion of global uniform exponential stability as considered in the previous section may be too strong for the point of view of applications.
Rather that requiring exponential decay of the solution for all possible signals $\sigma(\cdot)$, one may also be interested in guaranteeing
asymptotic convergence for some sufficiently large subset of possible signals. In this section, we review stability results for signals 
\begin{equation}
 \sigma\: [0,\infty) \to Q=\{0,1\}
\end{equation}
satisfying a persistent excitation (PE) condition, \ie, it holds
\begin{equation}
 \label{eq:Tmu}
\int_t^{t+T}\sigma(s)\,ds\geq \mu,\quad t \geq 0,
\end{equation}
for two positive real numbers $\mu \leq T$ independent of $t$. We briefly refer to $\sigma(\cdot)$ being a $T$-$\mu$ PE-signal, if \eqref{eq:Tmu} holds. 
Here we focus on system of the form
\begin{equation}\label{eq:sysPE}
\begin{aligned}
\frac{d}{dt} y(t) &= Ay(t)-\sigma(t) BB^*y(t),\quad t>0,\\
y(0)&=y_0,
\end{aligned}
\end{equation}
where we assume that $Y$ is a Hilbert space, $A$ is the generator of a strongly continuous semigroup of contractions $\SG(t)$ on $Y$,
$B$ is a bounded linear operator from a Hilbert space $U$ to $Y$ and $B^*$ denotes its adjoint.

The motivation for studying \eqref{eq:sysPE} is found in the fact that the linear feedback control $u(t)=-B^*y(t)$ is a common choice for stabilizing
a linear control system $\frac{d}{dt}y=Ay+Bu$ in the above setting going back to \cite{Haraux1989,Slemrod1974}. Hence, the dynamics of \eqref{eq:sysPE} consists of
switching between an uncontrolled evolution (for $\sigma(t)=0$) and a controlled one (for $\sigma(t)=1$). From this viewpoint, it seems natural to
impose conditions such as \eqref{eq:Tmu} on $\sigma(\cdot)$ guaranteeing a sufficient amount of action on the system if the uncontrolled system is
unstable, but asymptotic stability shall be achieved. For the analysis, we allow $\sigma(\cdot)$ taking also intermediate values 
$\sigma\: [0,\infty) \to Q=[0,1]$, noting that all results remain valid also for the case $Q~=~\{0,1\}$. 

If both the state and control space $Y$ and $U$ have finite dimension (and the pair of matricies $(A,B)$ is stabilizable), it is known that the switched system \eqref{eq:sysPE} is exponentially stable uniformly for the class of PE-signals $\sigma(\cdot)$ satisfying \eqref{eq:Tmu}; see \cite{ABJKKMPR,ChailletChitourLoriaSigalotti2008}. 
If $Y$ is infinite-dimensional, one can easily construct counter examples for such a result. We may, for instance, consider Example~\ref{ex:wave} for $N=1$ and $\Omega=(0,1)$ with $d=\chi_{(0,1)}$ for some proper subinterval $(a,b)$ of $\Omega$. Setting $b'=\frac{1+b}{2}$, $T=2$, and $\mu=1-b'$, the switching signal
\begin{equation}
 \sigma(t)=\sum_{k=0}^\infty \chi_{[2k-\mu,2k+\mu]}(t)
\end{equation}
satisfies \eqref{eq:Tmu} and the function
\begin{equation}
 v(t,x)=\sum_{k=0}^\infty\left(\chi_{[b'+2k,1+2k]}(x+t)-\chi_{[-1-2k,-b'-2k]}(x-t)\right)
\end{equation}
is a periodic, nonzero solution of \eqref{eq:wellengleichung-N} corresponding to $\sigma(\cdot)$ according to d'Alembert.
Similar examples can be found in \cite{HarauxMartinezVancostenoble2005,MartinezVancostenoble2002}. The given example also shows that also asymptotic stability (in a strong or weak sense) fails to hold.

Hence, in the general setting above, additional assumptions have to be imposed in order to guarantee the stability of the switched system and there are many previous studies devoted to conditions ensuring stability of second-order systems with time-varying parameters, mostly but not exclusively in the finite-dimensional setting, \eg, 
\cite{HarauxMartinezVancostenoble2005,Hatvani1996,PucciSerrin1996,Smith1961}. In \cite{HanteSigalottiTucsnak2012}, we have shown that exponential stability can be recovered if an appropriate observability inequality is satisfied. More precisely, the following result can be obtained.

\begin{theorem}\label{thm:PEexpstab} Under the assumption that 
\begin{equation}\label{eq:ineqneeded-intro}
 \int_0^\vartheta\sigma(t)\|B^* \SG(t)y_0\|_Y^2\,dt \geq c\|y_0\|_Y^2\quad\text{for all}~\mbox{$T$-$\mu$ PE-signal}~\sigma(\cdot)
\end{equation}
holds for some constants $\vartheta,c>0$ independent of $\sigma(\cdot)$, there exist two constants $M \geq 1$ and $\gamma>0$ independent of $\sigma(\cdot)$ such that the mild solution $y(\cdot)$ of system~\eqref{eq:sysPE} satisfies
\begin{equation}\label{eq:thm2est}
 \|y(t)\|_Y \leq M e^{-\gamma t}\|y_0\|_Y,\quad t\geq 0,
\end{equation}
for any initial data $y_0 \in Y$ and any $T$-$\mu$ PE-signal $\sigma(\cdot)$.
\end{theorem}
The proof of this result is based on deducing from \eqref{eq:ineqneeded-intro} a uniform decay for the solutions of \eqref{eq:sysPE} of the squared norm, chosen as Lyapunov function, on time-intervals of length $T$. The conclusion follows from considerations on the scalar-valued Lyapunov function as, for instance, earlier used in \cite{AeyelsPeuteman1998}. 

In the literature, the additional assumption \eqref{eq:ineqneeded-intro} is called \emph{generalized observability inequality} and is known to hold, for example, for the heat equation with boundary or locally distributed control \cite{Fattorini2005,MizelSeidman1997,Phung,Wang2004}. Using an eigenfunction expansion for the Dirichlet--Laplace operator, \cf~Example~\ref{ex:heat}, 
we have shown in \cite{HanteSigalottiTucsnak2012} that \eqref{eq:ineqneeded-intro} also holds for the $N$-dimensional wave equation from Example~\ref{ex:wave} if the damping satisfies $|d(x)| \geq d_0 >0$ for almost every $x \in \Omega$ and some constant $d_0$.

In \cite{HanteSigalottiTucsnak2012}, we have also shown that weak asymptotic stability can be recovered under weaker hypotheses using a compactness argument. 
\begin{theorem}\label{thm:PEweakstab}
Suppose that there exists $\vartheta>0$ such that for all $T$-$\mu$ PE-signals $\sigma(\cdot)$
\begin{equation}\label{eq:hypthm1}
 \int_0^\vartheta \sigma(t) \|B^* \SG(t)y_0\|^2_U\,dt = 0 \quad\Rightarrow\quad y_0 = 0.
\end{equation}
Then the mild solution $t \mapsto y(t)$ of system \eqref{eq:sysPE} converges weakly to the origin in $Y$ as $t\to\infty$ for any initial data $y_0 \in Y$ and any $T$-$\mu$ PE-signal $\sigma(\cdot)$.
\end{theorem}

Condition~\eqref{eq:hypthm1} is called a generalized unique continuation property, which corresponds for $\sigma(\cdot)\equiv 1$ to a property known to ensure approximate controllability; see, \eg, \cite{TucsnakWeiss2009}. In \cite{HanteSigalottiTucsnak2012} we show by an analyticity argument (Privalov's theorem) and standard unique continuation (Holmgren's theorem) that such generalized unique continuation properties hold for the Schr\"odinger equation in Example~\ref{ex:schoedinger} for damping $d \in L^\infty(\Omega)$ localized on an open nonempty subset $\omega\subset\Omega$ and satisfying
\begin{equation}
 |d(x)|\geq d_0 > 0,\quad x \in \omega~\text{almost everywhere}.
\end{equation}
Unless in the case $\omega=\Omega$, the question whether or not the stronger assumption~\eqref{eq:ineqneeded-intro} is satisfied in this example seems to be an open problem related to a conjecture of T.I. Seidman in \cite{Seidman1986}.

In the spirit of \cite{HarauxMartinezVancostenoble2005}, we also considered in \cite{HanteSigalottiTucsnak2012} excitations being rarefied in time and of variable duration rather than being satisfied on every time-window of prescribed length as in \eqref{eq:Tmu}. To this end we say that $\sigma(\cdot) \in L^\infty([0,T];[0,1])$ is of class $\mathcal{K}(A,B,T,c)$ if
\begin{equation}\label{eq:classK}
 \int_0^T \sigma(t)\|B^* \SG(t)y_0\|_U^2\,dt \geq c\|y_0\|_Y^2\quad \text{for all}~y_0 \in Y.
\end{equation}
The following result proved in \cite{HanteSigalottiTucsnak2012} then essentially states that strong stability is obtained when the total contribution of the excitations, suitably summed up, is infinite.
\begin{theorem}\label{thm:StrongStabAbstract}
Suppose that $(a_n,b_n)$, $n\in\NN$, is a sequence of disjoint intervals in $[0,\infty)$,
that $c_n$, $n\in\NN$, is a sequence of positive real numbers, and that
$\sigma(\cdot) \in L^\infty([0,\infty);[0,1])$ is such that its restriction $\sigma(a_n+\,\cdot~)|_{[0,b_n-a_n]}$ to
the interval $(a_n,b_n)$ is of class $\mathcal{K}(A,B,b_n-a_n,c_n)$ for all $n \in \NN$. Moreover, assume that 
\begin{equation}
\sup_{n\in\NN}(b_n-a_n) < \infty\quad\text{and}\quad\sum_{n=1}^\infty c_n = \infty. 
\end{equation}
Then the mild solution $t \mapsto y(t)$ of \eqref{eq:sysPE} satisfies $\|y(t)\|_Y \to 0$ as $t \to \infty$.
\end{theorem}

As an application, we may again consider the Schr\"odinger equation in Example~\ref{ex:schoedinger} for $N=1$, $\Omega=(0,1)$, and damping $d=\chi_\omega$, where $\omega$ is a nonempty
subinterval of $(0,1)$. It is well-known that for any interval $(a_n,b_n)$, $n\in\NN$, there exists a positive constant $c_n$ such that
\begin{equation}\label{eq:obsineq_n}
\int_{a_n}^{b_n} \int_{\omega} |\SG(t) y(x)|^2\,dx\,dt \geq c_n\|y\|_Y^2,~\quad y\in Y,
\end{equation}
that is, $\sigma(a_n+\cdot)|_{[0,b_n-a_n]}$ is of class $\mathcal{K}(A,B,b_n-a_n,c_n)$; see, for instance, \cite[Remark~6.5.4]{TucsnakWeiss2009}. Moreover, rewriting
\eqref{eq:obsineq_n} as
\begin{equation}\label{eq:obsineq2_n}
\int_{a_n}^{b_n} \int_{\omega} \left|\sum_{k \in \NN} \langle\phi_k,z\rangle_{L^2(0,1)}\phi_k(x) e^{i n^2 \pi^2 t} \right|^2\,dx\,dt \geq c_n\|y\|_Y^2,
\end{equation}
with $\phi_k(x)=\sqrt{2}\sin(n \pi x)$ we get from \cite[Proposition~7.1]{TenenbaumTucsnak2007} that $c_n$ can be taken satisfying
\begin{equation}\label{eq:ccost_sch1d}
c_n \geq C(b_n-a_n)^{-\frac{1}{2}}e^{-\frac{\pi}{2(b_n-a_n)}}
\end{equation}
for 
some positive constant $C$ independent of $n$. Hence, Theorem~\ref{thm:StrongStabAbstract} guarantees strong convergence for this example if
\begin{equation}
\sum_{n=1}^\infty (b_n-a_n)^{-\frac{1}{2}}e^{-\frac{\pi}{2(b_n-a_n)}}= \infty.
\end{equation}

\begin{remark} The results in the above example and, more generally, the methodology employed in this section, can be adapted to the case of some unbounded control operators and thus to boundary stabilization problems. We also note that the above conditions also provide a new method to obtain stability for intermittent damping in the finite-dimensional case
with integral excitations, where the sequence $c_n$ can be obtained using the Kalman rank, \ie, the minimal nonnegative integer $K$ such that 
\begin{equation}
 \text{rank}\,[B,AB,A^2B,\ldots,A^KB]=\dim(Y).
\end{equation}
Furthermore, we note that the results summarized above show that the sufficient condition for asymptotic stability of abstract second-order evolution equations with on/off damping considered in \cite{HarauxMartinezVancostenoble2005} is not a necessary condition. This question has been raised in \cite{FragnelliMugnai2008}.
Finally, we note that integral conditions in space, instead of in time as considered above, guaranteeing stabilizability of systems whose uncontrolled dynamics are given by a contraction semigroup have also been studied, \eg, \cite{Martinez1999,Nakao1996,Tcheugoue1998} for the wave equation and \cite{GuzmanTucsnak2003} for a plate equation. An interesting open question concerns
combined conditions in space-time to guarantee some sort of asymptotic stability.
\end{remark}

\subsection{Numerical results}\label{sec:numresultsstab} We briefly mention that the theoretical achievements presented so far can also be verified numerically. We demonstrate this using an example from \cite{AminHanteBayen2012} illustrating selected results from Section~\ref{sec:GUES} and point to applications.

\begin{figure}
\begin{tikzpicture}[]
\node[inner sep=0pt] (graph) at (0,0)
 {\includegraphics[trim=40 20 0 0,clip,width=.66\textwidth]{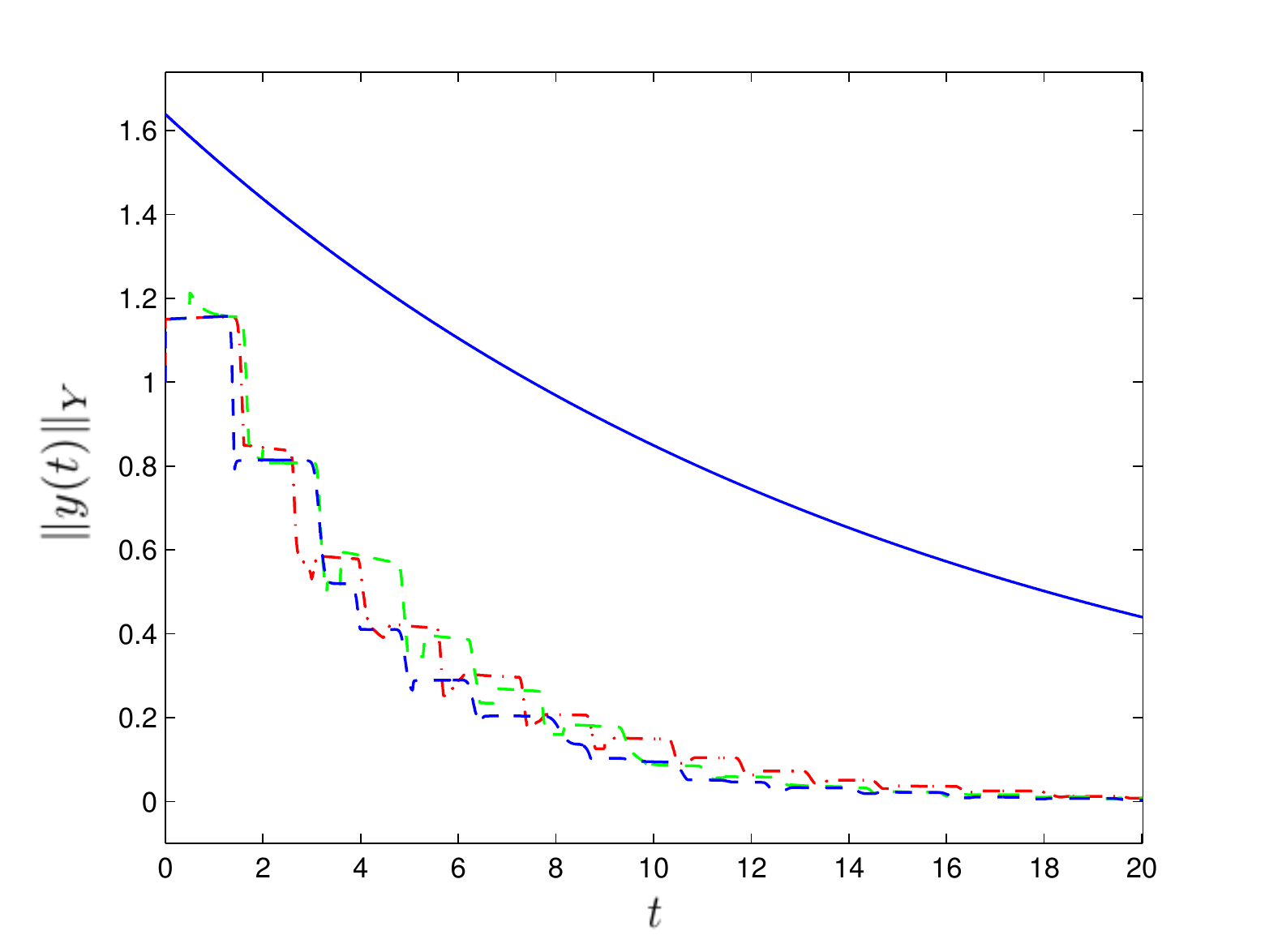}};
\node[inner sep=0pt,rotate=90] at (-4.5,0) {$\|y\|_\infty$};
\node[inner sep=0pt] at (-0.2,-3.4) {$t$}; 
\end{tikzpicture}
\caption{The exponential bound for $\|y(t)\|_\infty$ obtained in \cite{AminHanteBayen2012} (solid line) and the actual decay for three different switching signals (dashed lines) of solutions for Example~\ref{ex:transport} with the data from \eqref{eq:exsim}.}
\label{fig:sim}
\end{figure}

We consider Example~\ref{ex:transport} for a system of two equations and two modes, \ie, $n = 2$ and $M=2$.
The parameters and boundary data are specified as $[a,b]=[0,1]$, $f^j(y)=B^jy$, and
\begin{equation}\label{eq:exsim}
\begin{split}
        &\Lambda^{1}=\begin{bmatrix} -1.2 & 0\\ 0 & 1.8 \end{bmatrix}, 
\quad 	B^{1}=\begin{bmatrix} -0.005 & 0 \\ 0 & -0.005 \end{bmatrix},\\
	&\Lambda^{2}=\begin{bmatrix} -0.8 & 0\\ 0 & 1.4 \end{bmatrix},
\quad 	B^{2}=\begin{bmatrix} 0 & 0.005 \\ 0.005 & 0 \end{bmatrix},\\
&G_L^{1}=0.61,~G_R^1=1.15,~G_L^2=0.42,~G_R^2=1.21.
\end{split}
\end{equation} 
These matrices clearly satisfy \eqref{eq:stricthyperb}. Using that we have
\begin{equation}
\max_{j,j' \in Q} \rho\left(\begin{bmatrix} 0 & |G_R^{j'}|\\ |G_L^{j}| & 0 \end{bmatrix}\right) = 0.7381 < 1
\end{equation}
and $\|B^{1,2}\|_\infty=0.005$, we can conclude from Theorem~\ref{thm:GUEShyperbolic} and Remark~\ref{rem:GUEShyperbolicInhomo} that the switched system is globally uniformly exponentially stable. Figure~\ref{fig:sim} shows the predicted exponential bound from \cite{AminHanteBayen2012} and the actual decay of $\|y(t)\|_\infty$ for three different switching signals $\sigma(\cdot)$, where the solutions are approximations computed using the two-step Lax--Friedrichs finite difference scheme from \cite{Shampine2005}.

An application for the stability analysis of a cascade of water canals with open-close switching underflow slice gates using linearized Saint-Venant equations as a model is discussed in \cite{AminHanteBayen2008}.

\section{Optimal control of switching PDE-dynamical systems}\label{sec:optimalswitching} 
Consider an optimization problem of the form
\begin{equation}\label{eq:MIOCPorig}
\begin{array}{*1{>{\displaystyle\vphantom{\frac{d}{dt}}}l}}
 \text{minimize}~\varphi(y,u,\sigma) \; \text{subject to}\\
 \quad \frac{d}{dt} y(t)  =   A^{\sigma(t)}y(t)  + f^{\sigma(t)}(y(t),u(t)),\quad t \in (t_0,\, \tf),\\
 \quad y(0)= y_0,\\
 \quad g_k^\sigma(u,t) \leq 0~\text{for all}~t \in [t_0,\tf],~k=1,\ldots,m,\\
 \quad y \in C([t_0,\tf];Y),~u \in \Ut,~\sigma \in \Sigmat,
\end{array}
\end{equation}
where $\Ut$ is a Banach subspace of all measurable control functions $u\: [t_0,\tf] \to U$ for some Banach space $U$, $\Sigmat$ is a subspace of measurable integer control functions $\sigma\: [t_0,\tf] \to Q$, the dynamics for $y$ is a switched system of the form \eqref{eq:sysGeneral} augmented by an additional parameter $u$, and $\varphi$ is a cost function taking values in $\RR\cup\infty$. We assume that $Q$ is finite. Hence, problem \eqref{eq:MIOCPorig} is an optimal control problem with two controls $\sigma(\cdot)$ and $u(\cdot)$ which can be 
chosen independently in order to minimize the associated costs. The restriction that $\sigma(\cdot)$ takes values in a discrete set $Q$ reveals that \eqref{eq:MIOCPorig} can be regarded as an infinite-dimensional PDE-constraint mixed-integer type nonlinear optimization problem. Such problems cannot be solved with modern techniques from PDE-constraint optimization as they can be found, \eg, in \cite{HinzeEtAl2009}, because they typically rely on a (numerical) treatment of necessary optimality conditions in the spirit of the method of Lagrange multipliers, and this method is, in general, not available for mixed-integer type problems.

For the case that $A^{\sigma(\cdot)}$ is $\sigma(\cdot)$-invariant, we have therefore considered in \cite{HanteSager2013} and \cite{Hante2017} a relaxation of problem \eqref{eq:MIOCPorig}, which allows to explicitly construct integer-feasible controls being related to the original problem by the optimal value. The gap made by this approach can be estimated quite generally and can be shown to being arbitrarily small for particular problems. We review the essence of this method and the main technical arguments of the convergence analysis in Section~\ref{sec:relax}. Motivated by these results and again for the case that $A^{\sigma(\cdot)}$ is $\sigma(\cdot)$-invariant, we have further studied the optimal value function of the mixed-integer optimal control problem \eqref{eq:MIOCPorig}. In \cite{GugatHante2017}, we have shown that the optimal value is locally Lipschitz continuous as a function of problem parameters, somewhat naturally for perturbations of the initial data, and consistently with known results from mixed-integer linear programming for perturbations of the constraints on the controls in the case of linear quadratic problems under a Slater-type constraint qualification. The main results in these directions are reviewed in Section~\ref{sec:valuefunction}. For the case that $A^{\sigma(\cdot)}$ is $\sigma(\cdot)$-dependent and $u(\cdot)$ is fixed, we have studied the parametrization of $\sigma(\cdot)$ in problem~\eqref{eq:MIOCPorig} by switching times and mode sequences. With respect to these parameters, we have derived in \cite{RuefflerHante2016} differentiability properties of cost functions possibly involving switching costs subject to the switched PDE-dynamical system including discontinuous state resets at switching times and obtained gradient representation formulas based on an appropriate adjoint calculus. These results are reviewed in Section~\ref{sec:optimalityconditions} and can be used to solve mixed-integer optimal control problems of the form \eqref{eq:MIOCPorig} using gradient decent strategies.

\subsection{Outer convexification and relaxation for optimal switching control}\label{sec:relax}
Consider a mixed-integer optimal control problem of the form
\begin{equation}\label{eq:MIOCPSURorig}
\begin{aligned}
&\min~J=\phi(y(\tf)) \quad \text{s.\,t.}\\
&\quad \frac{d}{dt}y(t) = A y(t) + f^{\sigma(\cdot)}(y(t),u(t))~\text{on}~Y,~t \in (0,\tf),\\
&\quad y(0)=y_0, 
\end{aligned}
\end{equation}
where $A$ is assumed to be a generator of a strongly continuous semigroup of bounded linear operators $\SG(t)$ on a Banach space $Y$, $\tf>0$ is a fixed final time, and where, for all $j \in Q$, $f^j\: Y \times U \to Y$ and $\phi\: Y \to \RR$ are assumed to be continuous. The assumption that the cost function is of Mayer type is only for expository simplicity. 
More general Lagrange-type cost functions as~\eqref{eq:MIOCPorig} can either be treated using a variable transformation, or explicitly as done in \cite{HanteSager2013}.
Also, additional constraints as in \eqref{eq:MIOCPorig} can be handled. We discuss appropriate extensions in Remark~\ref{rem:SURconstraints}.

Associated with \eqref{eq:MIOCPSURorig}, consider the relaxed problem
\begin{equation}\label{eq:MIOCPSURrelaxed}
\begin{aligned}
&\min~\tilde J=\phi(y(\tf))\quad\text{s.\,t.}\\
&\quad \frac{d}{dt}y = A y + \sum_{j=1}^M \beta_j\,f^{j}(y,u)~\text{on}~Y,~t \in (0,\tf),\\
&\quad y(0)=y_0,\\
&\quad \sum_{j=1}^M \beta_j(t) = 1,~t \in (0,\tf),
\end{aligned}
\end{equation}
where $M=|Q|$ and the minimization is now with respect to the controls $u\: [0,\tf] \to U$ and $\beta=(\beta_1,\ldots,\beta_M)$ with each $\beta_j\: [0,\tf] \to [0,1]$. Observe that \eqref{eq:MIOCPSURrelaxed} is a PDE-constrained optimal control problem with ordinary constraints on the controls which can be assessed for a large class of operators $A$ with existing variational or operational methods; see, \eg, \cite{BensoussanDaPratoDelfourMitter2007,HinzeEtAl2009}. 

It is well-known that, under certain technical assumptions, the solution set of \eqref{eq:MIOCPSURorig} is dense in the solution set of \eqref{eq:MIOCPSURrelaxed} if $\sigma(\cdot)$, $u(\cdot)$, and $\beta(\cdot)$ are measurable control functions. This follows from recasting the problems \eqref{eq:MIOCPSURorig} and \eqref{eq:MIOCPSURrelaxed} using differential inclusions and applying the generalizations of the Filippov--Wa\v{z}ewski theorem proved in \cite{Frankowska1990} for $Y$ being a separable Banach space and in \cite{deBlasiPianigiani1999} for non-separable Banach spaces. While these results rely on powerful selection theorems, our main contribution in \cite{HanteSager2013,Hante2017} is a constructive proof of this result based on piecewise constant approximations of $\beta(\cdot)$ using binary values. Similar as in the case of ODE-constrained problems \cite{Sager2005}, this approach yields numerical methods leading to $\varepsilon$-optimal integer-feasible controls $\sigma(\cdot)$ for the original problem \eqref{eq:MIOCPSURorig}.

\begin{figure}
 \includegraphics[width=1.0\textwidth]{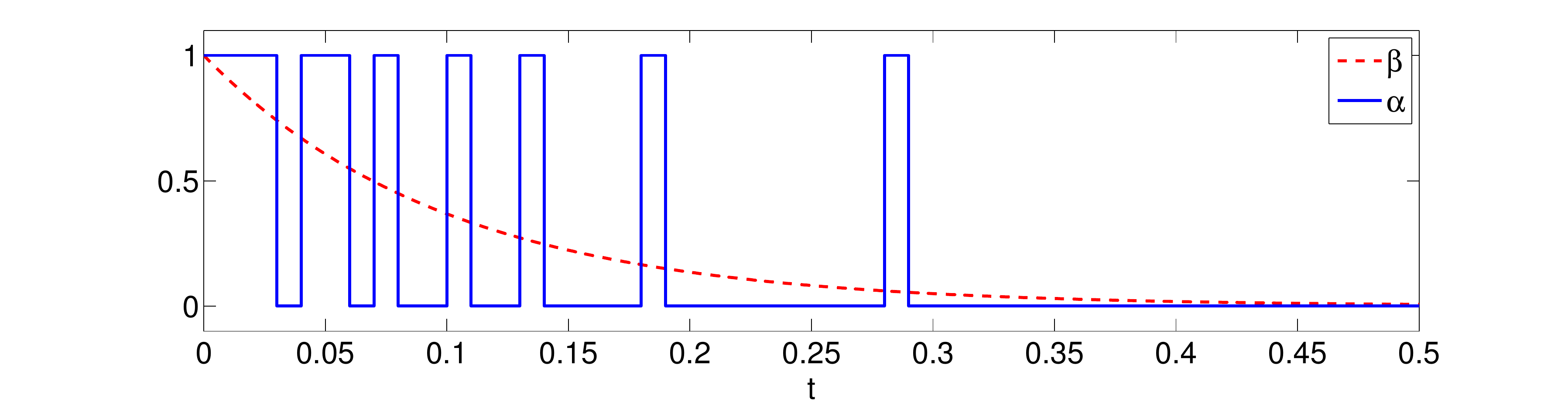}\vspace*{-2em}
 \caption{A scalar example of the binary function $\alpha$ obtained by the sum-up rounding strategy \eqref{eq:SUR1}--\eqref{eq:SUR2} for a singular arc $\beta(\cdot)$.}
\label{fig:SUR}
\end{figure}

To this end, let $u(\cdot)$ and $\beta(\cdot)$ be some measurable feasible controls for problem~\eqref{eq:MIOCPSURrelaxed}, let $0 = t_0 < t_1 < \cdots < t_n = \tf$ be some
given time instances, and consider defining a piecewise constant function $\alpha=(\alpha_1,\ldots,\alpha_M)\: [0,\tf] \to \{0,1\}^M$ by
\begin{equation}\label{eq:SUR1}
            \alpha_j(t) = p_{j,k}, \quad t\in [t_k,t_{k+1}),~j \in Q,~k=0,\ldots,n-1,
\end{equation}
where for all $j \in Q$, $k=0,\ldots,n-1$, $\hat{p}_{j,k}= \int_0^{t_{k+1}} \beta_j(\tau)\,d\tau - \sum_{l=0}^{k-1} p_{j,l} (t_{l+1}-t_{l})$, and
\begin{equation}\label{eq:SUR2}
\begin{aligned}
            &p_{j,k} = \begin{cases}1,~&\text{if}~\left(\hat{p}_{j,k} \geq \hat{p}_{l,k}~\forall~l \in Q \setminus\{j\}\right) \text{and}~\left(j<l~\forall~l \in Q \setminus\{j\}~\text{:}~\hat{p}_{j,k}=\hat{p}_{l,k}\right),\\
            0,~&\text{else}.\end{cases}
            \end{aligned}
\end{equation}
A scalar example of this sum-up rounding strategy is given in Figure~\ref{fig:SUR}. Then, defining a switching signal $\sigma(\cdot)$ on $[0,\tf]$ by
\begin{equation}\label{eq:SUR3}
\sigma(t)=\sum_{j=1}^M \alpha_j\,j,\quad t \in [0,\tf],
\end{equation}
the controls $\sigma(\cdot)$ and $u(\cdot)$ are feasible for the original problem \eqref{eq:MIOCPSURorig} and we can prove the following result concerning their quality as a suboptimal solution of \eqref{eq:MIOCPSURorig} in terms of the maximal grid size.

\begin{theorem}\label{thm:SURparabolic} Let $\bar{M}$ denote the growth bound of the semigroup $\SG(t)$, let $f^j$ be bounded with a constant $M_f$, Lipschitz continuous with a Lipschitz constant $L_f$ for all $j \in Q$ and let the cost function $\phi$ be Lipschitz continuous with a Lipschitz constant $L_\phi$. Furthermore, suppose that there exists a constant $C$ such that $\max_{j=1,\ldots,M} \|\frac{d}{ds}\SG(t-s)f^j(y(s),u(s))\| \leq C$ for almost every $t \in (0,\tf)$. Then it holds 
\begin{equation}
|\tilde J(u,\beta)-J(u,\sigma)| \leq \left(L_\phi(M_f + CT)e^{\bar{M} L_f \tf}\right) (M-1) \Delta t,
\end{equation} 
where $\Delta t=\max_{k=1,\ldots,n-1} (t_{k+1}-t_{k})$. 
\end{theorem}

If the technical hypotheses are satisfied, the above result applied to an optimal solution $u^*$ and $\beta^*$ of \eqref{eq:MIOCPSURrelaxed} particularly states that the optimal value of \eqref{eq:MIOCPSURrelaxed} can be reached arbitrarily close by the costs of \eqref{eq:MIOCPSURorig} associated to integer-feasible controls $u^*$ and $\sigma^*$ obtained via \eqref{eq:SUR1}--\eqref{eq:SUR3}.

The proof of Theorem~\ref{thm:SURparabolic} given in \cite{HanteSager2013} essentially uses that the binary approximation $\alpha$ of the function $\beta$ in \eqref{eq:SUR1}--\eqref{eq:SUR2} defining $\sigma(\cdot)$ via \eqref{eq:SUR3} satisfies $\sum_{j=1}^M \alpha_j(t) = 1$ for all $t \in [0,t_f]$, and that the integrated difference of $\alpha$ and $\beta$ satisfies the bound
\begin{equation}
  \max_{j=1,\ldots,M} \left| \int_0^t \alpha_j(\tau) - \beta_j(\tau)\,d\tau \right| \leq (M-1) \Delta t\quad~\text{for all}~t \in [0,t_f],
\end{equation} 
see \cite{SagerBockDiehl2011}. Further, it uses that under the assumed hypotheses, integration by parts applied to the variation of constants formula \eqref{eq:varconst} together with the Gronwall inequality yields a bound of the form
\begin{equation}\label{eq:solest0}
\begin{aligned}
 &\|y(t;u,\beta)-y(t;u,\sigma)\|_Y \\
 &\qquad \leq \left(L_\phi(M_f + CT)e^{\bar{M} L_f t}\right)\max_{j=1,\ldots,M} \left| \int_0^t \alpha_j(\tau) - \beta_j(\tau)\,d\tau \right|.
\end{aligned} 
\end{equation}

The assumption in Theorem~\ref{thm:SURparabolic} that $\frac{d}{ds}\SG(t-s)f^j(y(s),u(s))$ is bounded for almost every $t \in (0,\tf)$ can be verified for a large class of parabolic problems, \ie, when the semigroup $\SG$ generated by $A$ is analytic, and hence $s \mapsto \SG(t-s)$ is differentiable on $t-s>0$, and $f^j\: D(A) \to D(A)$. For hyperbolic problems, in general, such assumptions concerning the smoothness of the map $s \mapsto \SG(t-s)$ are not satisfied. In order to show that the above approach is also valid for a large class of hyperbolic problems, we
have studied the canonical system of Example~\ref{ex:transport} as above for the case that $\Lambda^{\sigma(\cdot)}$ (and hence $A^{\sigma(\cdot)}$) is $\sigma(\cdot)$-invariant, \ie, we consider the mixed-integer optimal control problem \eqref{eq:MIOCPSURorig} and the relaxation \eqref{eq:MIOCPSURrelaxed} for
\begin{equation}
\begin{aligned}
 &D(A)=\biggl\{\begin{pmatrix}y_I\\y_{II}\end{pmatrix} \in H^1(a,b;\RR^m) \times H^1(a,b;\RR^{(n-m)}) : \\
 &\qquad\qquad\qquad y_{II}(a)=G_L^j y_{I}(a),~y_{I}(b)=G_R^j y_{II}(b)\biggr\},\\ 
 &A y(x) = - \Lambda \frac{d}{dx}y,\quad x \in (a,b),~y \in D(A)
\end{aligned}
\end{equation}
with a continuously differentiable matrix function $\Lambda(x)=\mathrm{diag}\,(\lambda_1(x),\ldots,\lambda_n(x))$ satisfying $\lambda_1(\cdot)<\ldots<\lambda_{m}(\cdot)<0<\lambda_{m+1}(\cdot)<\ldots<\lambda_n(\cdot)$ on $[a,b]$, matrices $G_L$, $G_R$ of dimension ${(n-m)\times m}$ and ${m\times (n-m)}$, respectively, and a sufficiently smooth nonlinear functions $f^j$
for all $j \in Q$. We have shown in \cite{Hante2017}, that a result similar to the general one in Theorem~\ref{thm:SURparabolic} can be obtained in semi-classical Sobolev spaces introduced by \cite{Oberguggenberger1986} and being closely linked to the characteristic curves generated by the matrix $\Lambda$ given as solutions of the family of ODEs
\begin{equation}
\frac{d}{dt}s_i = \lambda_i(s_i),\quad s_i(\tau;\tau,\varsigma)=(\tau,\varsigma)\quad\text{for}~(\tau,\varsigma) \in \Omega_{\tf}=[0,\tf] \times [a,b]. 
\end{equation}

More precisely, for any $p,q,\mu \in \NN$ and any family of disjoint open sets $\Omega_m \subset \RR^p$, $m=1,\ldots,\mu$, let
\begin{equation}
\bigotimes_{m=1}^\mu W^{1,1}(\Omega_m;\RR^q)\quad~\text{and}\quad \bigotimes_{m=1}^\mu C^{0}(\Omega_m;\RR^q)
\end{equation}
denote the set of functions $h$ defined on the closure of $\bigcup_{m=1}^\mu \Omega_m$ with image in $\RR^q$ so that their restriction to $\Omega_m$
belongs to the classical Sobolev space $W^{1,1}(\Omega_m;\RR^q)$ or the Banach space $C^0(\Omega_m;\RR^q)$, respectively. 

Then, we make the following assumptions. Suppose that the initial data $y_0$ in problem \eqref{eq:MIOCPSURorig} is an element of the space $\bigotimes_{i}^\nu W^{1,1}((x_i,x_{i+1});\RR^n)$ for some partition $a=x_0<x_1<\ldots<x_{\nu-1}<x_\nu=b$ of the interval $[a,b]$.
For any piecewise smooth control $\beta\: [0,\infty) \to [0,1]^M$ with discontinuities at $\{\theta_i\}_{i=1}^\infty$, let $E^0(\beta)$ denote the union of all forward characteristic curves $s_i$ generated by $\Lambda$  (and their reflections at the boundaries) which lie in $\Omega_\infty = [0,\infty) \times [a,b]$, which emerge from the boundary points $\{(0,x_0),\ldots,(0,x_\nu\})$ and their intersection points with the sets $\{\theta_i\} \times [a,b]$ for all $i=1,\ldots,\infty$.
Further, for $k\in\NN$, let $E^k(\beta)$ be the union of $E^{k-1}(\beta)$ and all forward characteristic curves (and their boundary reflections) emerging from intersection 
points of characteristics that define $E^{k-1}$. Let $E_{\tf}(\beta)$ be the closure of all points in $\bigcup_{k=1}^\infty E^k(\beta) \cap \Omega_{\tf}$.

Supposing that
\begin{equation}\label{eq:HypNoDense}
 E_{\tf}(\beta) \cap \{t\} \times [a,b]~\text{is nowhere dense in}~\{t\} \times [a,b],\quad t\in (0,\tf],
\end{equation}
$E_{\tf}(\beta)$ is defined by finitely many discrete curves, which divide $\Omega_{\tf}$ up into finitely many simply connected open sets 
$\Omega_{\tf}^m$, $m=1,\ldots,\mu$, for some natural number $\mu$. We set
\begin{equation}
W^{1,1}_*(\Omega_{\tf} \setminus E_{\tf}(\beta)) := \bigotimes_{m=1}^\mu W^{1,1}(\Omega_{\tf}^m),\quad C^0_*(\Omega_{\tf} \setminus E_{\tf}(\beta)) := \bigotimes_{m=1}^\mu C^0(\Omega_{\tf}^m)
\end{equation}
and recall that the $L^1$-solution $y$ of the relaxed state equation in \eqref{eq:MIOCPSURrelaxed} can be obtained as a fixed point $y\: \Omega_T \to \RR^n$ of the integral transformation 
\begin{equation}
\psi(y)(\tau,\varsigma)=(\psi_1(y)(\tau,\varsigma),\ldots,\psi_n(y)(\tau,\varsigma)) 
\end{equation}
with $\psi_i(y)(\tau,\varsigma)$ defined recursively as
\begin{equation}\label{eq:psidef2}
 \psi_i(y)(\tau,\varsigma) = Y_i(\psi;\tau,\varsigma) + \sum_{j=1}^M \int_{t_i^*}^\tau \beta_j(t) f^j(y(t,s_i(t;\tau,\varsigma)),u(t))\,dt,
\end{equation}
where $t_i^*=t_i^*(\tau,\varsigma)$ denotes the intersection time of the curve $s_i(\cdot;\tau,\varsigma)$ with the boundary of $\Omega_T$
backward in time and $Y_i(\psi;\tau,\varsigma)$ being the $i$th component of the respective initial or boundary data at the intersection 
point. Moreover, under assumption \eqref{eq:HypNoDense} one can show that this solution satisfies 
\begin{equation}\label{eq:yreg}
 y \in W^{1,1}_*(\Omega_{\tf} \setminus E_{\tf}(\beta)) \cap C^0_*(\Omega_{\tf} \setminus E_{\tf}(\beta))
\end{equation}
for any piecewise smooth control $u\: [0,\infty) \to U$ and any piecewise smooth control $\beta\: [0,\infty) \to [0,1]^M$, see \cite{Oberguggenberger1986}.

We can then prove the following result concerning the relaxation gap.

\begin{theorem}\label{thm:SURhyperbolic} Let $\tf>0$ be sufficiently small (in the sense that \eqref{eq:HypNoDense} holds) and let $\phi$ be Lipschitz continuous. Further, let
$u\: [0,\tf] \to U$ and $\beta\: [0,\tf] \to [0,1]^M$ be piecewise smooth controls feasible for the relaxed problem \eqref{eq:MIOCPSURrelaxed} and let $\sigma(\cdot)$ be defined
by \eqref{eq:SUR1}--\eqref{eq:SUR3}. Then the controls $u(\cdot)$ and $\sigma(\cdot)$ are feasible for the original problem \eqref{eq:MIOCPSURorig} and there exists a constant $C>0$ (independent of $\Delta t$) such that 
\begin{equation}
|\tilde J(u,\beta)-J(u,\sigma)| \leq C \Delta t,
\end{equation} 
where $\Delta t=\max_{k=1,\ldots,n-1} (t_{k+1}-t_{k})$. 
\end{theorem}

The proof in \cite{Hante2017} uses a perturbation argument for the contraction mapping principle applied to the contractions 
$\psi=(\psi_1,\ldots,\psi_n)$ and $\varphi=(\varphi_1,\ldots,\varphi_n)$ given by \eqref{eq:psidef2} associated to the fixed control 
$(u,\beta)$ and the fixed control $(u,\alpha)$ with $\alpha$ taking the role of $\beta$, respectively, to show that it suffices to 
obtain existence of a constant~$\tilde{C}$ such that
\begin{equation}\label{eq:toshowNew}
\int_a^b |\psi_i(y)(t,\varsigma)-\varphi_i(y)(t,\varsigma)| d\varsigma \leq \tilde C \max_{j=1,\ldots,M} \left| \int_0^t \alpha_j(\tau) - \beta_j(\tau)\,d\tau \right|
\end{equation}
for all $i=1,\ldots,n$ and $t \in [0,\tf]$ in order to conclude the above result. The existence of $\tilde{C}$ is then obtained using integration by parts on the left-hand side of \eqref{eq:toshowNew} and bounding all appearing terms using the regularity property \eqref{eq:yreg}.

\begin{remark} The result in Theorem~\ref{thm:SURhyperbolic} extends also to inhomogeneous boundary conditions
\begin{equation}
 y_{II}(t,a)=G_L^{\sigma(t)} y_{I}(t,a) + d_L(t),~y_{I}(t,b)=G_R^{\sigma(t)} y_{II}(t,b)+d_R(t),\quad t \in (0,\infty),
\end{equation}
if there exist finitely many points $0=\tau_0<\tau_1<\ldots<\tau_{K-1}<\tau_K=\tf$ so that the boundary data satisfies 
$d_L \in \bigotimes_{i=1}^K W^{1,1}(\tau_i,\tau_{i+1};\RR^r)$ and 
$d_R \in \bigotimes_{i=1}^K W^{1,1}(\tau_i,\tau_{i+1};\RR^{(n-r)})$.
\end{remark}

\begin{remark}\label{rem:density} For the important case of hyperbolic systems of two variables assumption \eqref{eq:HypNoDense} is always satisfied \cite{Oberguggenberger1986}
and Theorem~\ref{thm:SURhyperbolic} holds for arbitrary $\tf>0$.
\end{remark}

\begin{remark}\label{rem:SURconstraints} In both Theorem~\ref{thm:SURparabolic} and \ref{thm:SURhyperbolic}, further constraints on $y$, $u$, and $\sigma(\cdot)$ can also be taken into account. Suppose that, for example, we wish to include a constraint of the form
\begin{equation}\label{eq:stateconstraint}
 g(y(t)) \leq 0,\quad t \in [0,t_f],
\end{equation}
for some Lipschitz continuous function $g\: Y \to \RR$ in the mixed-integer optimal control problem \eqref{eq:MIOCPSURorig}. Including this constraint also in 
\eqref{eq:MIOCPSURrelaxed} and assuming that a solution of this problem exists, the chain rule with the estimate \eqref{eq:solest0} then yields
\begin{equation}\label{eq:stateconstraintest}
|g(y(t;\beta) - g(y;\sigma)| \leq C \max_{j \in Q} \left| \int_0^t \alpha_j(\tau) - \beta_j(\tau)\,d\tau \right|,
\end{equation}
for some constant $C$ independent of $\sigma$. Hence, the violation of \eqref{eq:stateconstraint} can also be made arbitrarily small.
We may also wish to include combinatorial constraints, \eg, of the form
\begin{equation}\label{eq:cconstraints}
\#_{j^i \cvto j^k}(\sigma) \leq \bar{M}^{i,k},\quad i\in I,~k \in K,
\end{equation}
into the mixed-integer optimal control problem \eqref{eq:MIOCPSURorig}, 
where $\#_{j^i \cvto j^k}(\sigma)$ denotes the number of switches of the 
signal $\sigma\: [0,t_f] \to Q$ from value $j^i$ to value $j^k$, and 
$\bar{M}^{i,k}$ and $I,K \subset Q$ are given, nonnegative constants 
and index sets, respectively.  

We note that the sum-up rounding strategy \eqref{eq:SUR1}--\eqref{eq:SUR3}  
typically satisfies
\begin{equation*}
\#_{j^i \cvto j^k}(\sigma) \to +\infty\quad\text{for}~\Delta t \to 0 
\end{equation*}
for some $i,k \in Q$, so eventually violating \eqref{eq:cconstraints} for the 
desired accuracy.
Therefore, we may replace \eqref{eq:SUR1} and \eqref{eq:SUR2} by solving a mixed-integer program 
on the control discretization grid that minimizes 
\begin{equation}\label{eq:costminmax}
 \max_{j \in Q} \left| \int_0^t \alpha_j(\tau) - \beta_j(\tau)\,d\tau \right|
\end{equation}
with respect to $\alpha$ while taking into account the additional constraints \eqref{eq:cconstraints}.
This min-max problem can be written as a standard mixed-integer linear problem (MILP) using slack variables 
and can be solved using tailored branch-and-bound techniques \cite{SagerJungKirches2011}.
The resulting switching signal obtained via \eqref{eq:SUR3} then does, for $\Delta t \to 0$, in general 
not approach the optimal value of the relaxed problem \eqref{eq:MIOCPSURrelaxed}, but it yields a suboptimal
solution providing an upper bound for the optimal value of \eqref{eq:MIOCPSURorig} subject to \eqref{eq:cconstraints} 
via \eqref{eq:solest0} and the optimal value of the corresponding MILP.

Further details are discussed in \cite{HanteSager2013} and \cite{Hante2017}.
\end{remark}

\subsection{Lipschitz continuity of the value function}\label{sec:valuefunction}
Consider a mixed-integer optimal control problem with a parameter $\lambda$, 
\begin{equation}\label{eq:miocp}
\begin{array}{*1{>{\displaystyle\vphantom{\frac{d}{dt}}}l}}
 \text{minimize}~\varphi(\lambda,y,u,v) \; \text{subject to}\\
 \quad \frac{d}{dt} y(t)  =   A \, y(t)  + f^{\sigma}(t,y(t),u(t)),\quad t \in (t_0,\tf),\\
 \quad y(0)= y_0(\lambda),\\
 \quad g_k^{\sigma}(\lambda,u,t) \leq 0~\text{for all}~t \in [t_0,\tf],~k=1,\ldots,m,\\
 \quad y \in C([t_0,\tf];Y),~u \in \Ut,~\sigma \in \Sigmat,
\end{array}
\end{equation}
where $A$ is the generator of a strongly continuous semigroup on a Banach space~$Y$, and where $\Ut$ and $\Sigmat$ are
as defined in the beginning of Section~\ref{sec:optimalswitching}. Let $\nu(\lambda)\in \RR \cup \{\pm\infty\}$
denote the corresponding \emph{optimal value} given by
\begin{equation}\label{eq:defnu}
\begin{array}{*1{>{\displaystyle\vphantom{\frac{d}{dt}}}l}}
 \nu(\lambda)=\inf \bigl\{\varphi(\lambda,y,u,v) :\\
 \quad \frac{d}{dt} y(t)  =   A \, y(t)  + f^{\sigma}(t,\,y(t),u(t)),\quad t \in (t_0,\, \tf),\\
 \quad y(0)= y_0(\lambda),\\
 \quad g_k^{\sigma}(\lambda,u,t) \leq 0~\text{for all}~t \in [t_0,\tf],~k=1,\ldots,m,\\
 \quad y \in C([t_0,\tf];Y),~u \in \Ut,~\sigma \in \Sigmat\bigr\}.
\end{array}
\end{equation}
In \cite{GugatHante2017}, we have studied the dependency of $\nu(\lambda)$ on the parameter $\lambda$. The following adaption of a classical example 
shows that for switching control one can at most expect to have (local) Lipschitz continuity of the map $\lambda \to \nu(\lambda)$.
\begin{example}\label{ex:nonsmooth}
For some $\tf>0$ and $\lambda \in \RR$ consider the problem
\begin{equation}
\begin{array}{*1{>{\displaystyle\vphantom{\frac{d}{dt}}}l}}
 \text{minimize}~y(\tf)~\text{subject to}\\
 \quad \frac{d}{dt} y(t)=\sigma(t)\,y(t),~\text{for}~t \in (0,\tf),\quad y(0)=\lambda,\\
 \quad y(t) \in \RR,~\sigma(t) \in \{0,1\}~\text{for}~t \in (0,\tf).
\end{array}
\end{equation}
The optimal value function $\nu(\lambda)=\inf\{y(\tf;\lambda) : \sigma \in L^\infty(0,\tf;\{0,1\})\}$ can easily be seen to be
\begin{equation*}
 \nu(\lambda)=\begin{cases}e^{\tf}\lambda, &~\lambda<0,\\ \lambda, &~\lambda \geq 0, \end{cases}
\end{equation*}
which is Lipschitz continuous but not differentiable in $\lambda=0$. \qed
\end{example}
Our main results establish precisely this regularity for quite general mixed-integer optimal control problems when $\lambda$ acts solely on the initial data $y(0)$. 

\begin{assumption}\label{ass:generalLIP}
The function $\varphi\: \Lambda \times C([t_0,\tf];Y) \times \Ut \times \Sigmat \to \RR$ is continuous and, for every $\sigma \in \Sigmat$, the
functions $g^\sigma_1,\ldots,g^\sigma_m\: \Lambda \times \Ut \times [t_0,\tf] \to \RR$ are such that the set of admissible controls
\begin{equation}\label{defWT}
\begin{aligned}
 \Wt(\lambda):=\{&(u,\sigma) \in \Ut \times \Sigmat : \\
 &g_k^\sigma(\lambda,u,t) \leq 0,~k=1,\ldots,m,~t \in [t_0,\tf]\}. 
\end{aligned}
\end{equation}
is not empty for all $\lambda \in \Lambda$. 
The map $f^j\: [t_0,\tf] \times Y \times U \to Y$ is continuous for all $j \in Q$. Moreover, there exists a 
function $k \in L^1(t_0,\tf)$ such that for all feasible controls $(u,\sigma) \in \Wt(\lambda)$, for all $y_1,y_2 \in Y$ 
and for almost every $t \in (t_0,\tf)$
\begin{equation}
 |f^{\sigma(t)}(t,y_1,u(t))-f^{\sigma(t)}(t,y_2,u(t))| \leq k(t)|y_1 - y_2| 
\end{equation}
and
\begin{equation}
 |f^{\sigma(t)}(t,0,u(t))| \leq k(t).
\end{equation}
\end{assumption}

\begin{theorem}\label{thm:LipInitial} Under Assumption~\ref{ass:generalLIP}, suppose that
the constraint functions $g^\sigma_1,\ldots,g^\sigma_m$ are independent of $\lambda$. Let $\bar\lambda$ be some fixed parameter in $\Lambda$ and assume that for some
bounded neighborhood $B(\bar\lambda)$ of $\bar{\lambda}$ and some constant $L_0$
\begin{equation}\label{eq:y0Lip}
|y_0(\lambda_1)-y_0(\lambda_2)| \leq L_{0} \, |\lambda_1-  \lambda_2|,\,~\lambda_1, \lambda_2 \in B(\bar\lambda).
\end{equation}
Moreover, let $K=\sup_{\lambda \in B(\bar\lambda)}|y_0(\lambda)|$ and assume that for some constant $L_{\varphi}$
  \begin{equation}\label{eq:varphiLip}
   |\varphi(\lambda_1,y,u,\sigma)-\varphi(\lambda_2,\bar y,u,\sigma)| \leq L_{\varphi}(|y-\bar{y}|+|\lambda_1-\lambda_2|)
  \end{equation}
for all $u \in \Ut$ and $\sigma \in \Sigmat$ being feasible for problem \eqref{eq:miocp}, $y,\bar{y}$ such that $\max\{|y|,|\bar y|\} \leq C(\tf)(1+K)$ and
$\lambda_1$, $\lambda_2 \in B(\bar\lambda)$, where 
\begin{equation}
 C(t)=\gamma\exp\left(w_0 (t-t_0)+\gamma\int_{t_0}^{t}k(s)\,ds\right),
\end{equation}
for constants $\gamma \geq 0$ and $w_0 \geq 0$ such that $\|T(t)\| \leq \gamma \exp(w_0(t-t_0 ))$ for all $t \in [t_0,\tf]$.
Then there exists a constant $\hat L_\nu$ such that
\begin{equation}\label{eq:nuLipInitial}
 |\nu(\lambda_1)-\nu({\lambda_2})| \leq \hat L_\nu  |\lambda_1 - \lambda_2|,\quad~\lambda_1, \lambda_2 \in B(\bar\lambda).
\end{equation}
\end{theorem}

We have also obtained a similar result for $\lambda$ acting jointly on the initial data and the constraints for convex problems satisfying a Slater-type constraint qualification. 

\begin{assumption}\label{ass:Convex}
The mapping $(y,u) \mapsto f^{\sigma}(t,y,u)$ is linear and the mapping $(y,u) \mapsto \varphi(\lambda,y,u,\sigma)$ is convex.
Moreover, the function $\varphi$ is Lipschitz continuous with respect to $\lambda$ in the sense that
\begin{equation}
|\varphi(\lambda_1,y,u,\sigma) - \varphi(\lambda_2,y,u,\sigma)| \leq L_\varphi(|y|,|u|)|\lambda_1 - \lambda_2|
\end{equation}
with a continuous function $L_\varphi\: [0,\infty)^2 \rightarrow [0,\,\infty)$. 
For all $k=1,\ldots,m$, the mappings $u \mapsto g_k^{\sigma}(\lambda,u,t)$ are convex, the mappings $(u,t)\mapsto g_k^{\sigma}(\lambda,u,t)$ 
are continuous and the functions $g^{\sigma}_k$ are Lipschitz continuous with respect to $\lambda$ in the sense that for all 
$t \in [t_0,\tf]$
\begin{equation}
|g_k^{\sigma}(\lambda_1,u,t) - g_k^{\sigma}(\lambda_2,u,t)| \leq L_g(|u|)|\lambda_1 - \lambda_2|
\end{equation}
with a continuous function $L_g\: [0,\infty) \rightarrow [0,\,\infty)$.
\end{assumption}

\begin{theorem}\label{thm:JointLip} Under Assumptions~\ref{ass:generalLIP}~and~\ref{ass:Convex}, for any $\bar\lambda \in \Lambda$, and a bounded neighborhood $B(\bar\lambda) \subset \Lambda$ 
let $L_0,L_{\varphi}$ be constants such that \eqref{eq:y0Lip} and \eqref{eq:varphiLip} hold as in Theorem~\ref{thm:LipInitial}. 
Further, suppose that for some real numbers $\omega>0$ and some $\underline \alpha$ it holds that for all $\sigma \in \Sigmat$ there is a Slater point $\bar u_{\sigma} \in U$ 
such that for all $\lambda \in B(\bar \lambda)$ we have
\begin{equation}\label{eq:slaterpoint}
g_k^{\sigma}(\lambda,\bar u_{\sigma},t) \leq - \omega\quad\text{for all}~t \in [t_0,\tf],~k=1,\ldots,m,
\end{equation}
\begin{equation}
\label{30}
\sup_{\sigma \in \Sigmat}  \sup_{\lambda \in B(\bar \lambda)} \varphi(\lambda,y(\bar u_{\sigma},\sigma),\bar u_{\sigma},\sigma)<\infty,
\end{equation}
\begin{equation}
\label{31}
\nu(\lambda) \geq \underline \alpha
\end{equation}
and we have that the level sets
\begin{equation}
\label{coercive}
\hspace*{-1em}
\begin{aligned}
\bar S(y_0):=\bigcup_{\lambda_1,\lambda_2 \in B(\bar \lambda)}
\biggl\{& (u,\sigma) \in \Ut \times \Sigmat:\\
&~\varphi(\lambda_1,y(u,\sigma),u,\sigma) \leq \varphi(\lambda_1,y(\bar u_\sigma,\sigma),\bar u_\sigma,\sigma) + |\lambda_1 - \lambda_2|^2,\\
&~g_k^{\sigma}(\lambda_1,u(t),t) \leq 0,~t \in [t_0,\tf],~k=1,\ldots,m\biggr\}
\end{aligned}
\end{equation}
are uniformly bounded on $Y_0 = \{y_0(\lambda) : \lambda \in B(\bar\lambda)\}$.
Then there exists a constant~$L_\nu$ such that
\begin{equation}\label{eq:nuJointlyLip}
|\nu(\lambda_1)-\nu(\lambda_2)| \leq L_\nu |\lambda_1 - \lambda_2|\quad\text{for all}~\lambda_1,\,\lambda_2 \in B(\bar \lambda),
\end{equation}
where $\nu(\lambda)$ is the optimal value of \eqref{eq:miocp} as defined in \eqref{eq:defnu}.
\end{theorem}
The Slater-type constraint qualification in the latter result is in some sense a natural assumption and has a fully discrete counterpart in mixed-integer programming \cite{Williams1989}.
The proof relies on a strong duality result for parametric disjunctive programming.

In \cite{GugatHante2017}, we have shown that the assumptions of Theorem~\ref{thm:JointLip} are, \eg, satisfied for quadratic tracking-type optimal control problems
with a small uncontrollable disturbance $\varepsilon>0$ for $A$ defined as the Dirichlet--Laplace operator on a bounded domain $\Omega$ being controlled by switching the application 
of a lumped control $u \in [0,1+\varepsilon]$ between two non-overlapping control sub domains $\omega_1$ and $\omega_2$ with a constraint that restricts the values of $u$ to $[0,\varepsilon]$ 
for a dwell-time period of length~$\delta$ whenever a decision is taken to switch the control region. For this example, we have considered $\lambda$ being the joint 
perturbation of initial data, the disturbance $\varepsilon$, and the tracking target.

\subsection{Switching time and mode insertion gradients}\label{sec:optimalityconditions}
In this section, we consider problem \eqref{eq:MIOCPorig} for a moment with a continuous control $u$ being fixed, without control constraints (\ie, $m=0$) and study the case that the switching signal~$\sigma(\cdot)$ that is to be optimized acts on the generator and the nonlinear perturbation. More generally, we also include discontinuous state resets at switching times,
\ie, we study switching PDE-dynamics of the form
\begin{equation}
			\frac{d}{dt} y(t) = A^j y(t) + f^j(y(t)),\quad y = g^{j,j'} (y^-),
\end{equation}
whenever the mode $j \in \CM$ is held constant or whenever $j$ with associated state $y^-$ is switched to the new mode $j' \in \CM$ with new state $y$ at switching times $(\tau_k)_{k \in \NN} \subseteq [0,\tf]$, respectively. For our analysis, we consider a finite set of modes $Q$ and piecewise constant switching signals $\sigma(\cdot)$, which we parameterize by the sequences of switching times $(\tau_k)_k$ and modes $(j_k)_k$. The main result below concern gradient representation formulas for appropriate variations of these parameters based on solutions of adjoint problems, which are again a switched PDE-dynamical system. In this approach, we note that the costs $J$ may also include switching costs.

The hybrid semilinear evolutions are specified as follows: Given a fixed $N\in \NN_0$, a sequence of modes $j=(j_n)_{n=0,\ldots,N} \subseteq \CM$ and a monotonically increasing, but not necessarily strictly increasing sequence of switching times $\tau=(\tau_n)_{n=0,\ldots,N+1} \subseteq [0,\infty)$, we consider dynamics of the form
	\begin{alignat}{2}\label{eq:sys}
		\begin{aligned}
			\frac{d}{dt} y(t)	&= A^{j_n} y(t) + f^{j_n}(t,y(t)), \quad 	& n &\in \{0,\ldots,N\}, \, t \in (\tau_n,\tau_{n+1}),	\\
			y(\tau_n)	&= g^{j_{n-1},j_n}(y^-(\tau_n)),			& n &\in \{1,\ldots,N\},								\\
			y(\tau_0)	&= y_0.
		\end{aligned}
	\end{alignat}
A mapping $y\colon [\tau_0,\tau_{N+1}] \to Y$ is called a \emph{mild solution} of \eqref{eq:sys}, if, for all $n \in \{0,\ldots,N\}$, there are functions $y^n\colon [\tau_n,\tau_{n+1}] \to Y$ satisfying the following conditions:
\begin{enumerate}[label=(\roman*)]
\item The function $y^n$ is the only element of $C([\tau_n,\tau_{n+1}],Y)$ satisfying the variation of constants formula
	\begin{equation}
		y^n(t) = \SG^{j_n}(t-\tau_n)y^n_0 + \int_{\tau_n}^t \SG^{j_n}(t-s)f^{j_n}(s,y^n(s))\, ds	\quad~\text{for all}~t \in [\tau_n,\tau_{n+1}],
	\end{equation}
where
	\begin{equation}
		y^n_0 =
		\begin{cases}
			y_0, 								& \text{ if } n=0,	\\
			g^{j_{n-1},j_n}(y^{n-1}(\tau_n)),	& \text{ if } n \in \{1,\ldots,N\}.
		\end{cases}
	\end{equation}
\item If $\tau_n < \tau_{n+1}$ for some $n \in \{0,\ldots,N\}$, then $y|_{[\tau_n,\tau_{n+1})} \equiv y^n$.
\end{enumerate}
The map $y$ is called a \emph{classical solution} to \eqref{eq:sys}, if, furthermore, the following holds:
\begin{itemize}
\item[(iii)] If $\tau_n<\tau_{n+1}$ for some $n \in \{0,\ldots,N\}$, then $y^n \in C^1([\tau_n,\tau_{n+1}],Y)$.
\end{itemize}
We then define $y^-(\tau_n):=y^{n-1}(\tau_n)$ for all $n \in \{1,\ldots,N\}$. Depending on whether we wish to emphasize the dependence of a mild or classical solution $y$ to \eqref{eq:sys} on $(j,\tau)$ we use both the notations $y(\cdot)$ and $y(\cdot,j,\tau)$ equally in the following (however, keeping in mind not to confuse this with the value $y(\tau_k)=y(\tau_k,j,\tau)$ of $y$ at the time $t=\tau_k$).

Moreover, using induction over the number $N$ of switching points we have obtained the following well-posedness result.
\begin{lemma}\label{lem:existence} Under the general hypotheses introduced for system~\eqref{eq:sysGeneral} in Section~\ref{sec:intro}, there exists 
for any fixed number of switched $N \in \NN$, a unique maximal $\tf^\text{max}>0$ such that \eqref{eq:sys} has a unique mild solution on $[0,\tf^\text{max})$ for every sequence of modes $(j_n)_{n=0,\ldots,N} \subseteq \CM$ and every monotonically increasing sequence of switching times $(\tau_n)_{n=0,\ldots,N+1} \subseteq [0,\tf^\text{max})$. The maximal $\tf^\text{max}$ is lower semicontinuous as a function of the initial state $y_0 \in Y$. If, furthermore, $y_0 \in D(A^{j_0})$ and the mappings $g^{i,j}$ are continuously differentiable for all $i,j \in \CM$ with $i \neq j$, satisfying the inclusion $g^{i,j}(D(A^i)) \subseteq D(A^j)$, then the solution is classical.
\end{lemma}

Without loss of generality we set $\tau_0=0$ and, in regard of Lemma~\ref{lem:existence}, can add the assumption that $\tf \in (0,\tf^\text{max})$ is given with $T^\text{max}$ as in Lemma~\ref{lem:existence} and define the set of admissible switching times as
\begin{equation}
		\CT(0,\tf) = \{ \tau=(\tau_1,\ldots,\tau_N) \in \RR^N : 0 = \tau_0 \leq \tau_1 \leq \ldots \leq \tau_N \leq \tau_{N+1} = \tf \}.
\end{equation}
We may then consider a \emph{cost function} $J$ and the \emph{reduced cost function} $\Phi$ given by
\begin{align}
		\label{eq:cost}		J(\tau,y)		&= \int_0^{\tf} l(t,y(t))\, dt + \sum_{n=1}^N l^{j_{n-1},j_n}(\tau_n,y^-(\tau_n)),	\\
		\label{eq:redcost}	\Phi(j,\tau)	&= J(\tau,y(\cdot,j,\tau)),
\end{align}
and assume that $l \colon [0,\tf] \times Y \to \RR$ is continuous and continuously differentiable with respect to the second argument, and that $l^{m,n}\colon [0,\tf] \times Y \to \RR$ is continuously differentiable for every $m,n \in \CM$ with $m \neq n$. 

For a first result, we fix a sequence $j=(j_n)_{n=0,\ldots,N}$ of modes for the hybrid evolution \eqref{eq:sys} and address the subproblem of determining optimal switching times in order to minimize \eqref{eq:cost}. The problem can then be summarized as solving the following parametric optimization problem
	\begin{alignat}{3}\label{eq:min}
	 	\begin{aligned}
	 		\min\limits_{\tau} 	&		& J(\tau,y)\hspace{-0.25cm}		&																					\\
	 		\text{s.t.}			&		& \frac{d}{dt} y(t)		&= A^{j_n} y(t) + f^{j_n}(t,y(t)),\quad	& n &\in \{0,\ldots,N\},\, t \in (\tau_n,\tau_{n+1}), 	\\
	 							&		& y(\tau_n)			&= g^{j_{n-1},j_n}(y^-(\tau_n)),			& n &\in \{1,\ldots,N\},							\\
	 							&		& y(\tau_0)			&= y_0,																							\\
	 							&		& \tau				&\in \CT(0,\tf).
	 	\end{aligned}
	 \end{alignat}
Motivated by similar approaches for ODEs in \cite{DyerMcReynolds1970,EgWaAx06}, we have considered in \cite{RuefflerHante2016} the differentiability of $J$ with respect to admissible switching times $\tau \in \CT(0,\tf)$ and proved an adjoint equation based representation of the gradient $\frac{\partial \Phi}{\partial \tau}$. Analogous to the ODE case in \cite{EgWaAx06}, this leads to first-order optimality conditions and makes this subproblem \eqref{eq:min} accessible for gradient based optimization methods.

For a more detailed statement of the contribution, note that problem \eqref{eq:min} is equivalent to the minimization of the reduced cost function
\begin{equation}
		\Phi\colon \CT(0,\tf) \to \RR, \qquad \Phi(\tau)=J(\tau,y(\cdot,j,\tau))
\end{equation}
and we may conclude that a minimum exists using that $\Phi$ is continuous and $\CT(0,\tf) \subset \RR^N$ is compact. If $\Phi$ is even differentiable, we can ask for first-order optimality conditions. Formally applying the chain rule yields
\begin{equation}
		\frac{\partial \Phi}{\partial \tau} = \frac{\partial J}{\partial \tau} + \frac{\partial J}{\partial y}\frac{\partial y}{\partial \tau}.
\end{equation}
In order to evaluate the right-hand side by applying the chain rule, however, we would need to solve $N$ individual systems. Instead, we have investigated a computationally more efficient representation where we express the above derivative by means of the solution of \eqref{eq:sys} and the solution of the following \emph{adjoint problem} on the dual space $Y^*$: Find $p\colon [0,\tf] \to Y^*$ such that
\begin{alignat}{3}\label{eq:sysadj}
		\begin{aligned}
			\dot{p}(t)		&= -(A^{j_n})^* p(t) - [f^{j_n}_y(t,y(t))]^* p(t) +l_y(t,y(t)),\hspace{-3cm} 							&   & 							&   &						\\
							&																										& t &\in (\tau_n,\tau_{n+1}),	& n	&\in \{0,\ldots,N\},	\\
			p(\tau_n)		&= [g^{j_{n-1},j_n}_y(y^-(\tau_n))]^*p^+(\tau_n) -l^{j_{n-1},j_n}_y(\tau_n,y^-(\tau_n)),\hspace{-3cm}	&	&							&	&						\\
							&																										&	&							& n	&\in \{1,\ldots,N\},	\\
			p(\tf)			&= 0.																									&	&							&	&
		\end{aligned}
\end{alignat}

The adjoint equations can be motivated by a Lagrange formalism. Moreover, under our general hypotheses, we can see that the adjoint problem \eqref{eq:sysadj} has a unique mild solution if the forward problem \eqref{eq:sys} admits a classical solution. The first main result of \cite{RuefflerHante2016} is the following.

\begin{theorem}\label{theo:stgradient}
Assume $y$ is the unique classical solution of \eqref{eq:sys} and $p$ is the unique mild solutions of \eqref{eq:sysadj}. Then the reduced cost function $\Phi$ is continuously differentiable on $\CT(0,\tf)$ with respect to the $k$th switching time with
\begin{equation}
 \begin{aligned}
		\frac{\partial \Phi}{\partial \tau_k}(\tau)
		&= 				 l(\tau_k,y^-(\tau_k))-l(\tau_k,y(\tau_k))+l^{j_{k-1},j_k}_\tau(\tau_k,y^-(\tau_k))+\\
		&\phantom{{}=}	-\biggl\langle p^+(\tau_k), g^{j_{k-1},j_k}_z(z^-(\tau_k))\left(A^{j_{k-1}}z^-(\tau_k)+f^{j_{k-1}}(\tau_k,z^-(\tau_k))\right)+\\
		&\qquad\qquad\qquad\qquad -\left(A^{j_k} z(\tau_k)+f^{j_k}(\tau_k,z(\tau_k))\right) \biggr\rangle_\dual
 \end{aligned}
\end{equation}
for every $\tau \in \CT(0,\tf)$ and every $k \in \{1,\ldots,N\}$. 
\end{theorem}

The proof of the above result uses that for $\tau \in \CT(0,\tf)$ and $k \in \{1,\ldots,N\}$, the map $t \mapsto \left\langle p(t), y_k(t) \right\rangle_\dual $ defined for $t \in [\tau_k,\tau_{N+1}]$ is continuously differentiable on $(\tau_n,\tau_{n+1})$ for every $n \in \{k,\ldots,N\}$ with
\begin{equation}
		\frac{d}{dt} \left\langle p(t), y_k(t) \right\rangle_\dual = \left\langle l_y(t,y(t)), y_k(t) \right\rangle_\dual,
\end{equation}
where $y_k(t)$ is defined as the partial derivative $\frac{\partial y(t)}{\partial \tau_k}$ right-continuously extended on $[\tau_k,\tf]$ as the mild solution of the system
\begin{alignat}{2}\label{eq:sysdiff}
		\begin{aligned}
			\frac{d}{dt} y_k(t)	&= A^{j_n} y_k(t) + f^{j_n}_y(t,y(t)) y_k(t),								& t &\in (\tau_n,\tau_{n+1}),~n \in \{k,\ldots,N\},	\\
										y_k(\tau_n)		&= g^{j_{n-1},j_n}_y(y^-(\tau_n))y_k^-(\tau_n),								& n &\in \{k+1,\ldots,N\},						\\
			y_k(\tau_k)		&= g^{j_{k-1},j_k}_y(y^-(\tau_k))\left(A^{j_{k-1}}y^-(\tau_k)+f^{j_{k-1}}(\tau_k,y^-(\tau_k))\right)\hspace{-10cm}	&	& 	\\
							&\phantom{{}=} -\left(A^{j_k} y(\tau_k)+f^{j_k}(\tau_k,y(\tau_k))\right).	&	&
		\end{aligned}
\end{alignat}
A necessary condition can then be obtained by the classical necessary optimality conditions by Karush-Kuhn-Tucker for the constraint $\tau \in \CT(0,\tf)$.

\begin{theorem} \label{thm:noc} Under the assumptions of Theorem~\ref{theo:stgradient}, let $\tau$ be a local minimum of the reduced costs $\Phi$. Define
\begin{align*}
		a(\tau,n) &= \min \{ m \in \{0,\ldots,n\} \, | \, \tau_m=\tau_n \},	\\
		b(\tau,n) &= \max \{ m \in \{n,\ldots,N+1\} \, | \, \tau_m = \tau_n \}.
\end{align*}
Then it holds that
\begin{equation}
		\sum_{j=a(\tau,k)}^k \frac{\partial \Phi}{\partial \tau_j}(\tau) 	\leq 0 \qquad \text{and} \qquad
		\sum_{j=k}^{b(\tau,k)} \frac{\partial \Phi}{\partial \tau_j}(\tau) \geq 0
\end{equation}
for all $k \in \{1,\ldots,N\}$.
\end{theorem}

The results above apply for the case of a fixed sequence of modes. As a separate problem, we have in \cite{RuefflerHante2016} also studied the infinitesimal insertion of a new mode into a given sequence of modes for the hybrid evolution \eqref{eq:sys} and have again obtained an adjoint based representation for the sensitivity of the cost function \eqref{eq:cost} with respect to this perturbation. This concept has been introduced for ODEs in \cite{EgWaAx06} and makes the subproblem of determining optimal sequences of modes for the hybrid evolution \eqref{eq:sys} in order to minimize \eqref{eq:cost} again accessible for gradient based optimization methods. 

For a more detailed summary of the obtained result, assume that transition functions $g^{i,j}, g^{k,j}, g^{i,k}$ mapping between any modes $i,j,k \in \CM$ satisfy $g^{i,j}=g^{k,j} \circ g^{i,k}$, let $j=(j_n)_{n=0,\ldots,N} \subseteq \CM$ be a given sequence, $k \in \{0,\ldots,N\}$ be a fixed index and consider the insertion of the mode $\hat{\j} \in \CM$ at the time $\hat{\tau}=\tau_k$. Denote by 
\begin{align}\label{eq:expSwitching}
	\begin{aligned}
		j'		&= (j_1,\ldots,j_{k-1},\hat{\j},j_k,\ldots,j_N),	\\
		\tau'	&= (\tau_0,\ldots,\tau_k,\hat{\tau},\tau_{k+1},\ldots,\tau_{N+1})
	\end{aligned}
\end{align}
the expanded mode sequence and the switching time sequence, respectively, and denote by $y(\cdot,j',\tau')$ the solution of \eqref{eq:sys} with the additional mode, \ie, $y(\cdot,j',\tau')$ solves the expanded system
	\begin{alignat*}{2}
		\begin{aligned}
			\frac{d}{dt} y(t)		&= A^{j_n} y(t) + f^{j_n}(t,y(t)),		 			& n &\in \{0,\ldots,N\}\setminus\{k\}, \, t \in (\tau_n,\tau_{n+1}),		\\
			\frac{d}{dt} y(t)		&= A^{\hat{\j}} y(t) + f^{\hat{\j}}(t,y(t)),		& t &\in (\tau_k,\hat{\tau}),												\\
			\frac{d}{dt} y(t)		&= A^{j_{k+1}} y(t) + f^{j_{k+1}}(t,y(t)),\!\!\!	& t &\in (\hat{\tau},\tau_{k+1}),													\\
			y(\tau_n)		&= g^{j_{n-1},j_n}(y^-(\tau_n)),					& n &\in \{1,\ldots,N\}\setminus\{k+1\},									\\
			y(\hat{\tau})	&= g^{j_k,\hat{\j}}(y^-(\hat{\tau})),				&	&																		\\
			y(\tau_{k+1})	&= g^{\hat{\j},j_{k+1}}(y^-(\tau_{k+1})),			&	&																		\\
			y(\tau_0)		&= y_0.												&	&
		\end{aligned}
	\end{alignat*}
We distinguish between the adjoint solutions $p(\cdot,j,\tau)$ and $p(\cdot,j',\tau')$ in the same way. To indicate whether this expansion diminishes the cost function, we consider the \emph{mode insertion gradient}
	\begin{equation}\label{eq:migradient}
		\frac{\partial \Phi(\tau,j)}{\partial j_k} := \lim\limits_{\hat{\tau} \searrow \tau_k} \frac{J(\tau,y(\cdot,j',\tau'))-J(\tau,y(\cdot,j,\tau))}{\hat{\tau}-\tau_k}.
	\end{equation}
Under the assumptions of Theorem~\ref{theo:stgradient}, we have obtained the following result.
\begin{theorem}\label{theo:migradient} The mode insertion gradient \eqref{eq:migradient} is given by
	\begin{align}\label{eq:miformula}
		\begin{aligned}
			&\frac{\partial \Phi(\tau,j)}{\partial j_k}=  l(\tau_k,y(\tau_k,j,\tau))-l(\tau_k,y(\tau_k,j',\tau'))+(l^{k,\hat{\j}}_\tau(\tau_k,y(\tau_k,j,\tau))\\
			&\phantom{{}=}	+\biggl\langle p(\tau_k,j',\tau'), g^{j_k,\hat{\j}}_y(y(\tau_k,j,\tau))\left(A^{j_k}y(\tau_k,j,\tau)+f^{j_k}(\tau_k,y(\tau_k,j,\tau))\right)\\
			&\phantom{{}=}	-\left(A^{\hat{\j}} y(\tau_k,j',\tau')+f^{\hat{\j}}(\tau_k,y(\tau_k,j',\tau'))\right)\biggr\rangle_\dual.
		\end{aligned}
	\end{align}	
\end{theorem}

\begin{remark} The insertion of a new mode at some existent switching time $\tau_k$ is just for expository simplicity. The insertion of a new mode can actually be considered at any time $\hat{t}$ in $[0,\tf]$.  
\end{remark}

In \cite{RuefflerHante2016}, we consider as an example the energy-optimal switching from an unstable transport equation to an asymptotically stable diffusion equation and used Theorem~\ref{theo:stgradient} to verify that $\frac{\partial \Phi}{\partial \tau}\geq 0$. This yields the expected result that $\tau=0$ is a global minimum.
Theorem~\ref{theo:stgradient} and Theorem~\ref{theo:migradient} extend and, in a certain sense, unify the concept of switching-time optimization and mode insertion from switching ODE-dynamical systems in \cite{EgWaAx06} and ordinary DDE-dynamical systems in \cite{Ve2005,WuEtAl2006} to the abstract setting of nonlinearly perturbed strongly continuous semigroups. Unlike in most of the previously available work, the above results consider non-autonomous dynamics, state-resets at switching times, and include switching costs. Moreover, among switching of the nonlinear perturbation, the theory explicitly considers switching of the generators, which (in non-trivial cases) cannot be handled with the results available in the literature so far. This allows---under certain technical restrictions---the treatment of switching, \eg, the delay parameter of a DDE or switching the principle part of a PDE as in Example~\ref{ex:transport}. A conceptual algorithm with convergence analysis for an embedding of switching time and mode insertion gradients as obtained in Theorem~\ref{theo:stgradient} and Theorem~\ref{theo:migradient} using alternating directions, gradient projection and Armijo--Goldstein step size conditions can be found in \cite{AxelssonEtA2008}.

\subsection{Numerical results}\label{sec:numresultscontrol} Section~\ref{sec:relax} and Section~\ref{sec:optimalityconditions} are in principle closely related to numerical methods for optimal switching control problems. We have implemented these methods and tested them on several examples. To exemplify this, we summarize selected results from \cite{Hante2017} concerning the relaxation method for hyperbolic PDE-dynamical systems presented in Section~\ref{sec:relax} and again point to applications.

Consider the integer controlled nondiagonal semilinear hyperbolic system
\begin{equation}\label{eq:burgersrelax}
\begin{aligned}
 \eta_t + \xi_x = 0,\qquad \xi_t + a^2 \eta_x = -\kappa^{-1}(\xi-g^{\sigma(t)}(\eta)), \quad \sigma(t) \in Q=\{1,2\}
\end{aligned}
\end{equation}
subject to periodic boundary conditions $\eta(t,0)=\eta(t,L)$, $\xi(t,0)=\xi(t,L)$, where 
$g^1(\eta)=\frac12 \eta^2$, $g^2(\eta)=-\frac12 \eta^2$, $\kappa>0$, and $a^2$ such that $a^2-\eta^2\geq 0$. 
In characteristic variables $y_1 = \eta + a \xi$ and $y_2 = \eta - a \xi$, system \eqref{eq:burgersrelax} 
can be written as a diagonal system as in Example~\ref{ex:transport} with $\Lambda=\text{diag}(a,-a)$, $[a,b]=[0,L]$, 
and a nonlinear function $f^j(y)$. For sufficiently small $\kappa$, system \eqref{eq:burgersrelax} is an 
approximation of the control system 
\begin{equation}\label{eq:burgers}
 \eta_t \pm {\textstyle \frac12} \eta^2_x = 0,
\end{equation}
where the control consists of switching the sign in the flux function of the conservation law \eqref{eq:burgers},
in the sense that for fixed $j$, we have $\eta_t \pm {\textstyle \frac12} \eta^2_x = \kappa((a^2-\eta^2)\eta_x)_x$ 
up to second order in $\kappa$; see~\cite{Bianchini2001,JinXin1995}. 

As an optimization problem, we consider the minimization of a tracking-type cost functional
\begin{equation}\label{eq:costtrack}
 J(\eta)=\frac12\int_0^L |\eta(\tf,x)-\bar{\eta}(x)|^2\,dx
\end{equation}
for a given reference solution $\bar{\eta}$. Such flux switching control problems can be seen as an academic benchmark problem 
for traffic flow control. 

In \cite{Hante2017}, we have derived a mixed-integer nonlinear programming (MINLP) formulation for this problem based on a first-order explicit 
upwind-scheme applied to \eqref{eq:burgers} and a trapezoidal rule for \eqref{eq:costtrack}. Moreover, we have derived the relaxed problem
formulation \eqref{eq:MIOCPSURrelaxed}, which can equivalently be written as the following problem with a single control $\beta(t)$
\begin{equation}\label{eq:burgersrelaxrelaxocp}
\begin{aligned}
  \min_{\beta(t) \in [0,1]}~J(\eta) &\quad \text{subject to}\\
   &  \eta_t + \xi_x = 0,\\
   &  \xi_t + a^2 \eta_x = \kappa^{-1}(\beta(t)\eta^2- {\textstyle \frac12} \eta^2 - \xi),\\
   &  \eta(0,x)=\eta_0(x),~\xi(0,x)=0,\\
   & \eta(t,0)=\eta(t,L),~\xi(t,0)=\xi(t,L).
\end{aligned}
\end{equation}
For any piecewise $W^{1,1}$ initial data the hypothesis of Theorem~\ref{thm:SURhyperbolic} are satisfied, \cf Remark~\ref{rem:density}.
We may therefore conclude that, for our test data above, the relaxation gap can be made arbitrarily small.

We have solved the relaxed problem \eqref{eq:burgersrelaxrelaxocp} numerically using a gradient decent method, where
the derivative of the reduced cost function $\tilde J(\beta)=J(\eta(\beta))$ has been obtained as
\begin{equation}\label{eq:costtrackGrad}
 \tilde J'(\beta)=-\int_0^L \frac{q}{\kappa}[g(\eta,2)-g(\eta,1)]\,dx = \int_0^L \frac{q}{\kappa} \eta^2\,dx,
\end{equation}
where $(p,q)$ are the solutions of the following adjoint equations
\begin{equation}\label{eq:burgersrelaxadjoint}
\begin{aligned}
   &  -p_t - a^2 q_x = \kappa^{-1} q (2\beta(t)-1)\eta,\\
   &  -q_t - p_x = -\kappa^{-1}q,\\
   &  p(T,x)=-(\eta(T,x)-\bar{\eta}(x)),~q(T,x)=0,\\
   &  p(t,L)=p(t,0),~q(t,L)=q(t,0).
\end{aligned}
\end{equation}
For the numerical solutions of \eqref{eq:burgers}, \eqref{eq:burgersrelaxrelaxocp}, and \eqref{eq:burgersrelaxadjoint} we discretized $[0,L]$ by $N_x$ cells choosing a CFL-consistent time discretization with the CFL constant $\frac{1}{2}$ for $[0,\tf]$ and applied a first-order finite-volume in space and implicit Euler in time (IMEX) splitting scheme from \cite{BandaHerty2009,JinXin1995} for the discretization of the PDE in \eqref{eq:burgersrelaxrelaxocp} and of \eqref{eq:burgersrelaxadjoint}.

As test data, we have taken in \cite{Hante2017} the parameters specified 
as $L=2\pi$, $T=3$, $\xi(0,x)=2\chi_{(\frac{L}{4},\frac{3L}{4})}(x)$, $x \in [0,L]$, and $\eta(0,x)=\bar{\eta}(x)=1-\sin(x)$, $x \in [0,L]$.
In this case, GAMS/BONMIN applied as a heuristic in order to solve the MINLP reports as the best found solution $\lambda(t)=1-\chi_{(\bar{t},T]}(t)$ 
with $\bar{t}\approx\frac94$. The computation time, however, grows exponentially with increasing $N_x$, and becomes excessive already for $N_x=70$, 
see Table~\ref{table:miocphyperb}. GAMS/BARON providing globally optimal solution confirms this solution for $N_x=10$, but runs out of 
CPU-time (12 hours) for $N_x=20$ and larger.

\begin{figure}
 \includegraphics[width=\textwidth]{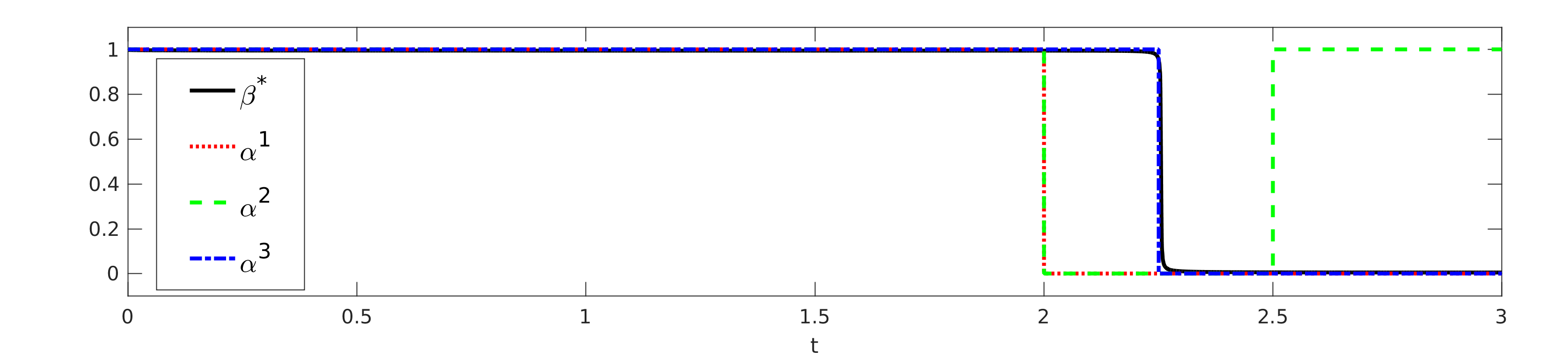} \vspace*{-2em}
\caption{Numerical results for a test problem from \cite{Hante2017}. The figure shows the optimal relaxed control $\beta$ and the 
integer control approximations $\alpha^1$, $\alpha^2$, and $\alpha^3$ corresponding to $\Delta t=1.0$, $\Delta t=0.5$, and
$\Delta t=0.25$ for $N_x=1000$. The relative integer optimality gaps $\gamma(\sigma^k)=\gamma(\alpha^k)=(J^*-J(\sigma^k))/J^*$ 
are $\gamma(\alpha^1)=4.4$, $\gamma(\alpha^2)=3.62$ and $\gamma(\alpha^3)=0$ in accordance with Theorem~\ref{thm:SURhyperbolic}.}
\label{fig:miocphyperb}
\end{figure}

The integer control found by relaxation and rounding then decreases with $\Delta t$ in accordance with the 
rate predicted in Theorem~\ref{thm:SURhyperbolic}. The computation time grows linearly, see Table~\ref{table:miocphyperb}, and the best found 
solution for $N_x=1000$ coincides with the MINLP-solution for $\Delta t=0.25$ and smaller on equidistant control grids, \cf~Figure~\ref{fig:miocphyperb}.
Results for an application of this method to the optimal switching of variable speed limit signs for a traffic flow benchmark problem from \cite{Hegyi2003}
using Jin–Xin relaxed Lighthill–Whitham–Richards (LWR) network model are available in \cite{Hante2016}.

\begin{table}[t]\small
\begin{center}
{\sc GAMS/BONMIN}\\
\begin{tabular}{cccccccc}
\toprule
$N_x$&	10 & 20 & 30 & 40 & 50 & 60 & 70 \\	
\midrule
\#vars&	187&756& 1\,705& 3\,034&4\,743&6\,832&9\,301\\
$J^*$&	 0.280&  0.211&  0.125&0.121&0.095&0.095&--\\
CPU-time (s)&0.74&8.44&163.17&734.95&632.75&2\,232.57&$>259\,200.00$\\
\bottomrule
\end{tabular}
\\[1em]
{\sc Outer Convexification/Relaxation}\\
\begin{tabular}{cccccc}
\toprule
$N_x$&	100 & 200 & 300 & 400 & 500 \\
\midrule
\#vars&	475\,165 & 1\,905\,151 & 4\,290\,539& 7\,630\,326& 11\,924\,913 \\
$J^*$&	 0.1952&  0.1109&  0.0875& 0.0785& 0.0748\\
CPU-time (s)&21.36&28.93&64.14&74.10& 117.71\\
\midrule
$N_x$&	600 & 700 & 800 & 900 & 1\,000 \\	
\midrule
\#vars&	17\,175\,501 & 23\,379\,888 & 30\,539\,075 & 38\,654\,863 &47\,723\,850\\
$J^*$&	0.0733& 0.0729 & 0.0731 & 0.0731 &0.0742\\
CPU-time (s)& 139.81 & 166.74 & 234.79 & 304.11 &318.86\\
\bottomrule
\end{tabular}
\\[1em]
\end{center}
\caption{Numerical results for a test problem from \cite{Hante2017} comparing a MINLP formulation with 
a numerical method supported by the theory summarized in Section~\ref{sec:relax}.}
\label{table:miocphyperb}
\end{table}

\section{Conclusion}\label{sec:conclusion}
The summarized results all regard abstract switched systems in Banach or Hilbert spaces. They address stability properties that hold uniformly with 
respect to switching and optimal switching control for linear parabolic as well as hyperbolic PDE-dynamical systems including linear and in some cases 
bounded nonlinear perturbations. 

The main contributions concern the theory of common Lyapunov functions and global uniform exponential stability, a theory of generalized observability inequalities 
for asymptotic stability of intermittently damped dissipative systems, a theory concerning the relaxation gap, the regularity of
the optimal value function, as well as an adjoint calculus for perturbations of switching times and mode insertions in 
optimal switching control. The contributions have been illustrated on examples including the heat, the Schrödinger's, the
wave, and the transport equation and selected results have also been verified numerically by case studies. 

Several individual open problems have been remarked. More general open problems concerning the stability of switched systems are necessary and sufficient conditions 
that are possibly less general but more convenient to be checked. Such conditions may combine for instance Lyapunov techniques with 
conditions on intermittencies. A similarly general open question concerning optimal switching control problems are the use of relaxation 
techniques that verifiably yield exact or $\varepsilon$-optimal solutions for problems involving switching costs.

\section*{Acknowledgments}
The author acknowledges the support from the DFG under the grant CRC/Transregio~154 enabling the work in \cite{Hante2017}, \cite{GugatHante2017} and \cite{RuefflerHante2016}, 
the support from the French program ARPEGE under the ANR grant ArHyCo contract number ANR-2008 SEGI 004 01-30011459 leading to the work in \cite{HanteSigalotti2011} and \cite{HanteSigalottiTucsnak2012}, and the support from the Mathematics Center Heidelberg (MATCH) facilitating the work in \cite{HanteSager2013}.



\end{document}